\documentclass[openany]{book}

\usepackage{import}

\usepackage[latin1]{inputenc}
\usepackage[english]{babel}

\usepackage{amsthm}
\usepackage{amsmath}
\usepackage{amssymb}
\usepackage{amsfonts}
\usepackage{mathtools}

\usepackage{graphicx}
\usepackage{float}
\usepackage{tikz-cd}
\usetikzlibrary{babel}

\usepackage[a4paper]{geometry}
\usepackage{fancyhdr}
\usepackage{scrextend}
\usepackage{afterpage}

\newcommand\blankpage{%
    \null
    \thispagestyle{empty}%
    \newpage}

\usepackage{mathrsfs}
\usepackage{dsfont}

\usepackage{enumerate}
\usepackage{enumitem}

\usepackage{nameref}
\usepackage{hyperref}
\usepackage[english]{cleveref}

\usepackage{thmtools}

\usepackage{titlesec}

\usepackage[backend=biber, style=alphabetic, citestyle=alphabetic, maxalphanames=4, maxnames=10]{biblatex} 
\usepackage{tocvsec2}

\hypersetup{
    colorlinks=true,
    linkcolor=blue,
    citecolor=blue
}

\declaretheoremstyle[
	qed=$\triangledown$
]{definitions}
\declaretheoremstyle[
	qed=$\blacktriangledown$
]{examples}
\declaretheoremstyle[
  	headfont=\normalfont\itshape,
	bodyfont=\normalfont,
  	postheadspace=1em,
]{remarks}

\declaretheorem[name=Definition,
	style=definitions, 
	numberwithin= chapter,
	refname={definition, definitions},
	Refname={Definition, Definitions},
	]{definition}

\declaretheorem[name=Example,
	style=examples, 
	sibling=definition,
	refname={example, examples},
	Refname={Example, Examples},
	]{example}
\declaretheorem[name=Theorem, 
	sibling=definition,
	refname={theorem, theorems},
	Refname={Theorem, Theorems}
	]{theorem}
\declaretheorem[name=Proposition, 
	sibling=definition,
	refname={proposition, propositions},
	Refname={Proposition, Propositions}
	]{proposition}
\declaretheorem[name=Lemma, 
	sibling=definition,
	refname={lemma, lemmas},
	Refname={Lemma, Lemmas}
	]{lemma}
\declaretheorem[name=Corollary, 
	sibling=definition,
	refname={corollary, corollaries},
	Refname={Corollary, Corolaries}
	]{corollary}
	
\newtheorem*{theorem*}{Theorem}
	
\newenvironment{Proof}
	{\begin{addmargin}[3em]{0em}
	\begin{proof}}
	{\end{proof}
	\end{addmargin}
	\vspace{6pt}}
	
\newenvironment{SProof}
	{\begin{addmargin}[3em]{0em}
	\begin{proof}[Sketch of the proof]}
	{\end{proof}
	\end{addmargin}
	\vspace{6pt}}

\titleformat{\part}[display]
	{\normalfont}
	{\Large\thepart\filcenter}
	{10pt}
	{\Huge\bfseries\filcenter}

\titleformat{\chapter}[display]
	{\normalfont}
	{\Large\thechapter\filcenter}
	{10pt}
	{\huge\bfseries\filcenter}
\titlespacing*{\chapter}
	{0pt} 
	{4cm} 
	{6cm} 

\titleformat{\section}[block]
  {
  \normalfont}
  {\Large\thesection.}
  {10pt}
  {\Large\bfseries}
  [\titlerule]
\titlespacing*{\section}
	{0pt} 
	{20pt} 
	{10pt} 
 
\titleformat{\subsection}[block]
  {
  \normalfont}
  {\large\thesubsection.}
  {10pt}
  {\large\bfseries}
\titlespacing*{\subsection}
	{0pt} 
	{10pt} 
	{10pt} 

\newcommand{\R}{\mathbb{R}} 
\newcommand{\C}{\mathbb{C}} 

\newcommand{\tr}{\operatorname{tr}} 

\newcommand{\im}{\operatorname{im}} 

\newcommand{\norm}[1]{\left \Vert #1 \right \Vert} 
\newcommand{\prodesc}[2]{\left \langle #1,#2 \right \rangle} 
\newcommand{\bracket}[2]{\left \{ #1,#2 \right \}} 

\newcommand{\Diff}{\operatorname{Diff}}

\newcommand{\Ad}{\operatorname{Ad}}
\newcommand{\ad}{\operatorname{ad}}

\newcommand{\GL}{\operatorname{GL}}
\newcommand{\SO}{\operatorname{SO}}

\newcommand{\pr}{\mathrm{pr}}

\newcommand{\gfrak}{\mathfrak{g}}
\newcommand{\tfrak}{\mathfrak{t}}
\newcommand{\hfrak}{\mathfrak{h}}
\newcommand{\Xfrak}{\mathfrak{X}}

\newcommand{\zfrak}{\mathfrak{z}}
\newcommand{\gl}{\mathfrak{gl}}
\newcommand{\so}{\mathfrak{so}}
\newcommand{\Lcal}{\mathcal{L}}
\newcommand{\Fscr}{\mathscr{F}}

\newcommand{\muhat}{\hat{\mu}}

\newcommand{\Cinf}{C^\infty}
\newcommand{\Mimpl}{M_{\mathrm{impl}}}
\newcommand{\muimpl}{\mu_{\mathrm{impl}}}
\newcommand{\muhatimpl}{\hat{\mu}_{\mathrm{impl}}}

\newcommand{\Xline}{\overline{X}}
\newcommand{\Yline}{\overline{Y}}
\newcommand{\Vline}{\overline{V}}

\newcommand{\Piline}{\overline{\Pi}}
\newcommand{\dbar}{/\!\!/\!}
\newcommand{\gammadot}{\dot{\gamma}}
\newcommand{\gdot}{\dot{g}}
\newcommand{\betadot}{\dot{\beta}}
\newcommand{\Gammadot}{\dot{\Gamma}}

\newcommand{\derev}[3]{\left.\frac{d #1}{d #2} \right\vert_{#2 = #3}} 
\newcommand{\derpar}[2]{\frac{\partial #1}{\partial #2}} 
\newcommand{\der}[2]{\frac{d #1}{d #2}} 

\begin{document} 

\maxtocdepth{section} 

\frontmatter
\frontmatter
\begin{titlepage}
\begin{center}
	\includegraphics[scale=.25]{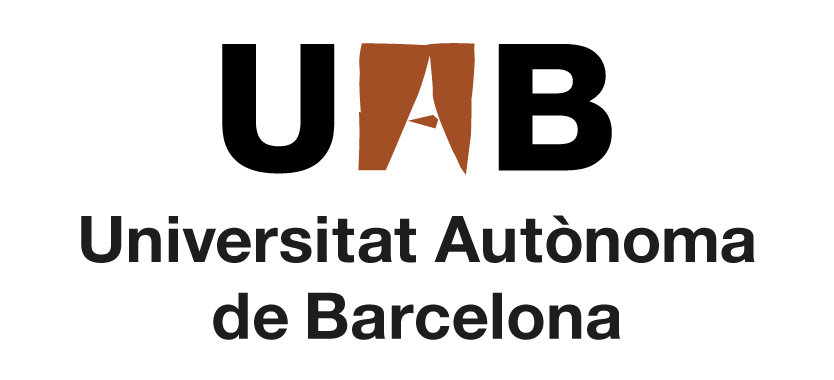} \par
	\vspace{5\baselineskip}
	\rule{\textwidth}{1.5pt}
	{\LARGE \textbf{From Symplectic to Poisson.}\par {\Large \textbf{A Study of Reduction and a Proposal Towards Implosion.}}\par}
	\rule{\textwidth}{1.5pt}
	\vspace{2\baselineskip}
	{\Large Bachelor's thesis \par}
	
	{\LARGE \[
	\begin{tikzcd}
  	(\mu^{-1}(\alpha),i^*\omega=\pi^*\omega_\alpha) \arrow[r,"i"] \arrow[d,"\pi",swap] & (M,\omega) \\
	(M\dbar_\alpha G,\omega_\alpha)
	\end{tikzcd}
	\]}
	
	\vspace{3\baselineskip}
	{\large Universitat Aut\`{o}noma de Barcelona \par \textit{Department of Mathematics} \par}
	\vspace{\baselineskip}
	{\large July 2021 \par}
	\vspace{5\baselineskip}
	\vspace{2\baselineskip}
	\textbf{Jaime Pedregal Pastor} \par
	\textit{Supervised by} Dr. Roberto Rubio \par
\end{center}
\end{titlepage}

\blankpage

\newpage
\thispagestyle{plain}
\begin{center}
    \Large
    \textbf{From Symplectic to Poisson.} \par
    \large \textbf{A Study of Reduction and a Proposal Towards Implosion}
    
    \vspace{0.4cm}
    Jaime Pedregal Pastor
    
    \vspace{0.9cm}
    \textbf{Abstract}
\end{center}

The imploded cross-section of a symplectic manifold is a stratified space allowing for an abelianization of its symplectic reduction. After recalling symplectic and Poisson reduction and reviewing the basics of symplectic implosion, we prove a cross-section theorem for Poisson manifolds, generalizing the Guillemin-Sternberg theorem for symplectic manifolds, which constitutes a first step towards Poisson implosion. On our way, we find and fix a mistake in the proof of Guillemin-Sternberg's theorem, and we identify Poisson transversals as the right analogue to symplectic submanifolds in this context.

\newpage
\thispagestyle{plain}
\begin{center}
    \large
    \textbf{Acknowledgements}
\end{center}

My deepest gratitude goes to my supervisor, Dr. Roberto Rubio. It has been a pleasure to work with him. I really appreciate his dedication and his time throughout the last year. He was always willing to meet me and discuss my advances on the project, and these discussions were very enriching. The effort he has made in revising the many preliminary versions of this dissertation cannot be sufficiently highlighted. I am specially thankful for his willingness to descend to details when it was needed, despite the arduous computations this may demand sometimes; for a student just starting to do research somewhat independently, this is immensely valuable. In addition, he has kindly advised me on my academic future and has shown true interest in it, something which clearly goes beyond the call of duty. He also opened my eyes to certain fields in geometry which were totally unknown to me and which now constitute one of my principal sources of inspiration. Thank you for all of that, Roberto.

Of course, I am also indebted to my family and my friends, who patiently supported me during the writing of this thesis. Many thanks to them also.

\tableofcontents


\mainmatter
\chapter*{Introduction}
\addcontentsline{toc}{chapter}{Introduction}
\markboth{Introduction}{Introduction}

Symplectic and Poisson geometries were originally conceived to provide a proper geometrical framework for the description of classical mechanics. Nevertheless, they soon became fields of independent mathematical interest and they have been extensively studied since then.

It is well known in physics that symmetries in a mechanical system give rise to conserved quantities in the dynamics of the system, a result commonly referred to as Noether's theorem. The existence of these conserved quantities, and the remaining symmetry on them, can be used to reduce the number of degrees of freedom of the mechanical system. Mathematically, symmetries are encoded as structure-preserving smooth actions of a Lie group $G$ on a symplectic or Poisson manifold $M$. Of special interest are the Hamiltonian actions, which come equipped with a smooth function $\mu:M\to\gfrak^*$, where $\gfrak$ is the Lie algebra of $G$, satisfying a certain dynamical condition. The map $\mu$ is called a moment map, because the classical linear and angular momentums can be recovered as moment maps for particular actions. In this case, $M$ can be ``reduced'' using $\mu$.  The precise meaning of this reduction will be given in \Cref{sect:symplectic-reduction}, but, loosely speaking, the level sets of the map $\mu$ can be quotiented by a subgroup of $G$ such that the quotient space is a symplectic or Poisson manifold whose structure is related to the one in $M$ by means of the projection. The precise mathematical formulation of reduction was firstly given by J. Marsden and A. Weinstein for the symplectic case \cite{marsden-weinstein} and by Marsden and T. Ratiu for the Poisson case \cite{marsden-ratiu}.

Besides the obvious application to solving mechanical systems with symmetry, the reduction process also offers a natural way of obtaining symplectic and Poisson structures on certain objects whose structures otherwise seem somewhat accidental. We explain in detail, for instance, the case of the Kirillov-Kostant-Souriau symplectic form on the coadjoint orbits of any Lie group, or the case of the Lie-Poisson structure on the dual of any Lie algebra.

V. Guillemin, L. Jeffrey and R. Sjamaar introduced in \cite{guillemin-jeffrey-sjamaar} the concept of the imploded cross-section of a symplectic manifold $M$ carrying a Hamiltonian action of a Lie group $G$. The imploded cross-section is a topological space $\Mimpl$ such that the reduction of $M$ by $G$ is the same as the reduction of $\Mimpl$ by a maximal torus $T$ of $G$. Hence, it somehow ``abelianizes'' the $G$-action on $M$. This does not come without a cost, though: the space $\Mimpl$ is singular, meaning that it is not a symplectic manifold, but rather a disjoint collection of symplectic manifolds---a stratification, to be precise. 

One of the main tools one uses to prove that the imploded space $\Mimpl$ is indeed a collection of symplectic manifolds is the Guillemin-Sternberg cross-section theorem, {\cite[Thm. 26.7]{guillemin-sternberg}}, which states that the preimage by $\mu$ of a submanifold transverse to a coadjoint orbit is locally a symplectic submanifold of $M$. The proof they provide in the book, nevertheless, contains a mistake at the end, where they implicitly assume the moment map to be a submersion, which is not true in general (see \Cref{ex:moment-maps-are-not-submersions}). The proof can be easily corrected to hold for the general case, and this is one of our contributions (see the proof of \Cref{thm:guillemin-sternberg-cross-section}).

Our main contribution, however, is a generalization of the Guillemin-Sternberg cross-section theorem to the Poisson setting. When trying to define a ``Poisson imploded cross-section,'' mimicking the symplectic construction, one would hope that it also decomposes as a disjoint union of Poisson manifolds. The preimage of the transverse submanifold is now expected to be some kind of submanifold of $M$ inheriting a Poisson structure from that of $M$, such that in the symplectic case it corresponds to a symplectic submanifold. We find that the correct notion is that of a Poisson transversal (\Cref{sect:poisson-transversals}). If $M$ is now a Poisson manifold carrying a Hamiltonian $G$-action with moment map $\mu$, and if we write $G\cdot\lambda$ for the coadjoint orbit of $\lambda\in\gfrak^*$, we prove the following result:

\newtheorem*{thm-poisson-cross-section}{Theorem \ref{thm:poisson-cross-section}}

\begin{thm-poisson-cross-section}[Poisson Cross-section]
Let $\lambda\in\gfrak^*$ and $Z\subseteq \gfrak^*$ a submanifold perfectly transverse to $G\cdot \lambda$ at $\lambda$, that is, such that $T_\lambda Z \oplus T_\lambda (G\cdot \lambda)=\gfrak^*$. Then, if $p\in M$ with $\mu(p)=\lambda$, there is a neighborhood $U$ of $p$ such that $\mu^{-1}(Z)\cap U$ is a Poisson transversal of $M$.
\end{thm-poisson-cross-section}

After proving the theorem, we describe the pathway to follow in order to complete the definition of a Poisson implosion with properties analogous to that of the symplectic implosion. It principally consists in the generalization to the Poisson setting of two results in symplectic geometry---the refinement of the cross-section theorem provided by E. Lerman et al. (\Cref{thm:lmtw-cross-section}) and the generalization of the reduction theorem of Sjamaar and Lerman (\Cref{thm:sjamaar-lerman})---and the use of Poisson reduction by stages.
\bigskip

The structure of the thesis is as follows. In \Cref{chap:moment-maps} we define symplectic and Poisson structures on manifolds and the notion of the moment map arising from a Hamiltonian group action, together with its properties and some examples. In \Cref{chap:reduction} we study reduction in both the symplectic and the Poisson contexts, exploring their dynamical counterpart and some interesting examples. Finally, in \Cref{chap:implosion} we first present the symplectic implosion construction, showing its partition into symplectic manifolds and making precise how it abelianizes the action, and then we introduce the Poisson transversals and prove the generalization of the Guillemin-Sternberg theorem. We have added an appendix as a reference for our notation and some relevant notions.
\bigskip

The Einstein summation convention is consistently used in the text: indices that appear repeated, once as a superindex and once as a subindex, are implicitly summed over all their possible values; e.g., the expression $p_idq^i$ is to be read as $\sum_{i=1}^n p_idq^i$. All manifolds are assumed to be smooth, meaning $\Cinf$. Given a local chart $(U,(x^i))$ on a manifold, we write $dx^i$ and $\partial_i$ for the coordinate 1-forms and vector fields, respectively.

\chapter{Moment maps}\label{chap:moment-maps}

We begin by reviewing the language we will be using in the rest of the dissertation: that of symplectic and Poisson manifolds, smooth Lie group actions and moment maps. We offer explicit constructions of moment maps in some special cases and several examples.


\section{Symplectic Geometry} 

Let $M$ be a smooth manifold and $\omega\in\Omega^2(M)$. We define the kernel of $\omega$ as the set
\[
\ker\omega:=\{v\in TM: i_v\omega=0\}\subseteq TM
\]
and say that $\omega$ is \textbf{nondegenerate} if $\ker\omega=0$. This is equivalent to the condition that the bundle morphism $\omega:TM\to T^*M$ given by $\omega(v):=i_v\omega$ is in fact a bundle isomorphism. We use the same symbol for the 2-form and the associated bundle morphism because it will be clear from the context which of the two we are referring to each time.

In local coordinates $(x^i)$ we can write
\[
\omega=\frac{1}{2}\omega_{ij}dx^i\wedge dx^j, \quad\text{with $\omega_{ij}:=\omega(\partial_i,\partial_j)$.}
\]
Nondegeneracy of $\omega$ amounts to the matrix $(\omega_{ij})$ being nonsingular at each point.

\begin{definition}
A \textbf{symplectic manifold} is a pair $(M,\omega)$, where $M$ is a smooth manifold and $\omega$ is a closed nondegenerate 2-form on $M$. If $(M',\omega')$ is another symplectic manifold, we say that a diffeomorphism $F:M\to M'$ is a \textbf{symplectomorphism} if $F^*\omega'=\omega$. If such a symplectomorphism exists we say that $M$ and $M'$ are symplectomorphic.
\end{definition}

If $(M,\omega)$ is a symplectic manifold, then necessarily the dimension of $M$ is even, since every odd-dimensional skew-symmetric matrix is singular.

Let $p\in M$ and let $W\subseteq T_pM$ be a vector subspace. We define its orthogonal space with respect to $\omega$ as
\[
W^\perp:=\{v\in T_pM: \omega(v,w)=0,\text{ for all $w\in W$}\}.
\]

Depending on the relation between $W$ and its orthogonal, we can consider different types of submanifolds on a symplectic manifold.

\begin{definition}
Let $(M,\omega)$ be a symplectic manifold. Let $p\in M$ and let $W\subseteq T_pM$ be a vector subspace. Then $W$ is said to be
\begin{enumerate}
	\item symplectic if $W\cap W^\perp=0$,
	\item isotropic if $W\subseteq W^\perp$,
	\item coisotropic if $W^\perp\subseteq W$, and
	\item Lagrangian if $W=W^\perp$.
\end{enumerate}

If $N\subseteq M$ is a submanifold, we say that $N$ is symplectic (resp. isotropic, coisotropic, Lagrangian) if $T_pN$ is a symplectic (resp. isotropic, coisotropic, Lagrangian) subspace of $T_pM$ for every $p\in N$.
\end{definition}

To study the dynamics on symplectic manifolds, we define two special kinds of vector fields on $M$.

\begin{definition}
Let $(M,\omega)$ be a symplectic manifold. We say that $X\in\Xfrak(M)$ is a \textbf{symplectic vector field} if $\Lcal_X\omega=0$. For any $f\in\Cinf(M)$, we define the associated \textbf{Hamiltonian vector field} as the unique vector field $X_f$ such that $i_{X_f}\omega=df$, i.e., $X_f=\omega^{-1}(df)$. We say that $X\in\Xfrak(M)$ is Hamiltonian if there is some $f\in\Cinf(M)$ such that $X=X_f$, in which case we say that $f$ is an associated Hamiltonian or energy function (it is not unique, although any two differ by a locally constant function).
\end{definition}

By Cartan's magic formula, $X\in\Xfrak(M)$ is a symplectic vector field if and only if $i_X\omega$ is closed, and it is Hamiltonian if in addition $i_X\omega$ is exact. We say that a flow on $M$ is Hamiltonian if the induced vector field is Hamiltonian.

\begin{example} \label{SympGeom:example:r2n}
Consider $\R^{2n}$ with Cartesian coordinates $(q^1,\dots,q^n,p_1,\dots,p_n,)$. Then the 2-form $\omega_0=dq^i\wedge dp_i$ is a symplectic form on $\R^{2n}$, called the standard symplectic form. A vector field $X\in\Xfrak(\R^{2n})$ is Hamiltonian with energy $H$ if and only if
\[
X=\derpar{H}{p_i}\derpar{}{q^i}-\derpar{H}{q^i}\derpar{}{p_i}.
\]
Its integral curves are precisely the solutions to Hamilton's equations as generally written in classical mechanics, see {\cite[\S40]{landau-lifshitz}}.
\end{example}

In fact, the Darboux theorem ensures that every symplectic form can be locally written in the standard form just described. For a proof see, for instance, {\cite[Sect. 3.2]{mcduff-salamon}} or {\cite[Thm. 22.13]{lee}}.

\begin{theorem}[Darboux]
Let $(M,\omega)$ be a symplectic manifold of dimension $2n$. For every $p\in M$ there is a chart $(U, (q^i,p_i)_{i=1}^n)$, with $p\in U$, such that $\omega=dq^i\wedge dp_i$ on $U$.
\end{theorem}

Another important example is the following.

\begin{example}
Let $M$ be any smooth manifold. Its cotangent bundle $T^*M$ has a natural symplectic structure. Let $\pi: T^*M\to M$ be the projection and define the tautological 1-form $\theta\in\Omega^1(T^*M)$ as
\[
\theta_\alpha(v):=\alpha(\pi_*v), \text{ for $\alpha\in T^*M$ and $v\in TT^*M$.}
\]
Then $\omega=-d\theta$ is a symplectic form on $T^*M$, called the canonical symplectic form on the cotangent bundle. To see that it is symplectic, let $(U,(x^i))$ be a chart on $M$ and let $(\pi^{-1}(U),(q^i,p_i))$ be the associated natural chart on $T^*M$, given by $q^i=x^i\circ \pi$ and $p_i(\alpha)=\alpha(\partial/\partial x^i)$. Then it is easy to see that $\theta=p_idq^i$, and hence $\omega=dq^i\wedge dp_i$.
\end{example}

This example allows us to extend the notion of classical mechanics on $\R^{2n}$, as in \Cref{SympGeom:example:r2n}, to classical mechanics on an arbitrary manifold $M$. If $M$ is the model of the space of degrees of freedom of a mechanical system, one can define a smooth function $H$ on the phase space $T^*M$, called the Hamiltonian or energy function, governing the dynamics of the system through its associated Hamiltonian vector field $X_H$. This makes sense because $T^*M$ is canonically endowed with a symplectic structure. The motions of the system are precisely the integral curves of $X_H$. 

This notion can be generalized even further, without asking for the phase space to be the cotangent bundle of the space of degrees of freedom: we call the triple $(M,\omega,H)$ a \textbf{Hamiltonian system} if $(M,\omega)$ is a symplectic manifold and $H\in\Cinf(M)$. Further information on symplectic geometry can be found in \cite{mcduff-salamon,cannasdasilva}.


\section{Group Actions}

Let $M$ be a smooth manifold and $G$ a Lie group acting smoothly on $M$. Let $\gfrak$ be the Lie algebra of $G$ and $\Phi$ the action of $G$ on $M$. We write $\Phi_g$ for the left translation by $g\in G$ on $M$ and $\Phi^p$ for the orbit map of $p\in M$. For any $\xi\in\gfrak$ we define the associated \textbf{infinitesimal generator} or fundamental vector field as the vector field on $M$ defined by
\begin{equation} \label{eq:def-infinitesimal-generator}
\xi_M(p):=\Phi^p_*(\xi)=\derev{}{t}{0}\Phi_{\exp(t\xi)}(p).
\end{equation}
We then define $\gfrak_M(p):=\{\xi_M(p):\xi\in\gfrak\}$, the set of infinitesimal generators at $p$.

There are two natural actions of $G$ on its Lie algebra $\gfrak$ and its dual $\gfrak^*$: the \textbf{adjoint} and \textbf{coadjoint} actions, respectively. To define them, let $C_g$ stand for conjugation by $g\in G$ on $G$, i.e., $C_g(h):=ghg^{-1}$. Then the adjoint action is given by $\Ad(g)\xi:=C_{g*}\xi$, for $\xi\in\gfrak$, and the coadjoint action by $\Ad^*(g)\alpha:=\Ad(g^{-1})^*\alpha$, for $\alpha\in\gfrak^*$. It is easy to check that these indeed define smooth actions of $G$ on $\gfrak$ and $\gfrak^*$. It is a well-known fact that the pushforward of the adjoint action $\Ad:G\to \GL(\gfrak)$ at the identity is
\[
\begin{array}{rrcl}
\ad: & \gfrak & \longrightarrow & \gl(\gfrak) \\
& \xi & \longmapsto & [\xi,\cdot]
\end{array}.
\]

The following proposition gives some properties concerning infinitesimal generators and these actions, which will be useful in computations.

\begin{proposition} \label{Group-Actions:inf-gen:properties}
For any $\xi,\eta\in\gfrak$ and $g\in G$ we have that
\begin{enumerate}
	\item $(\Ad(g)\xi)_M=\Phi_{g*}\xi_M$, and \label{Group-Actions:inf-gen:properties:1}
	\item $[\xi_M,\eta_M]=-[\xi,\eta]_M$. \label{Group-Actions:inf-gen:properties:2}
\end{enumerate}
\end{proposition}

\begin{Proof}
To see \cref{Group-Actions:inf-gen:properties:1}, we compute for any $p\in M$: 
\begin{align*}
(\Ad(g)\xi)_M(p) &=\derev{}{t}{0}\Phi_{\exp(tC_{g*}\xi)}(p)=\derev{}{t}{0}\Phi_{g\exp(t\xi)g^{-1}}(p) \\
&= \Phi_{g*}\left( \derev{}{t}{0}\Phi_{\exp(t\xi)}(\Phi_{g^{-1}}(p))\right)=\Phi_{g*}(\xi_M(\Phi_g^{-1}(p)))\\
&=(\Phi_{g*}\xi_M)(p).
\end{align*}
To see \cref{Group-Actions:inf-gen:properties:2}, using that $\Ad_*=\ad$ and \cref{Group-Actions:inf-gen:properties:1}, we obtain
\begin{align*}
[\xi_M,\eta_M](p)&= \Lcal_{\xi_M}\eta_M(p)=\derev{}{t}{0}(\Phi_{\exp(-t\xi)*}\eta_M)(p) \\
&= \derev{}{t}{0}(\Ad(\exp(-t\xi))\eta)_M(p)=\Phi^p_*\left( \derev{}{t}{0}\Ad(\exp(-t\xi))\eta\right) \\
&=(-\ad\xi(\eta))_M(p)=-[\xi,\eta]_M(p). & \qedhere
\end{align*}
\end{Proof}

The following proposition gives information about the structure of the orbits of the $G$-action. Notice that for any $p\in M$, the isotropy group $G_p$ is a closed subgroup of $G$ because $G_p=(\Phi^p)^{-1}(p)$. Hence the quotient $G/G_p$ is a smooth manifold such that the projection $\pi:G\to G/G_p$ is a smooth submersion (see {\cite[Cor. 4.1.21]{abraham-marsden}} or {\cite[Thm. 21.17]{lee}}).

\begin{proposition} \label{Group-Actions:injective-immersion}
Let $p\in M$, then the map
\[
\begin{array}{rrcl}
\Theta^p: &G/G_p & \longrightarrow & M \\
&gG_p & \longmapsto & \Phi^p(g)
\end{array}
\]
is an injective immersion. In particular, $\Theta^p(G/G_p)=G\cdot p$ is an immersed submanifold of $M$ such that $\Theta^p$ is a diffeomorphism onto its image.
\end{proposition}

\begin{Proof}
We have that $\Theta^p\circ\pi=\Phi^p$. This implies that $\Theta^p$ is smooth because $\Phi^p$ is smooth. Indeed, as an easy consequence of the local normal form for submersions, $\pi$ has smooth local sections about each point in $G/G_p$, and this immediately gives the smoothness for $\Theta^p$.

The map $\Theta^p$ is injective because $\Phi_g(p)=\Phi_h(p)$ implies $h^{-1}g\in G_p$. To see that it is an immersion, write $[g]=gG_p$ for simplicity. Since $\Theta^p\circ\pi=\Phi^p$, then $v=\pi_*w\in T_{[g]}(G/G_p)$ is in $\ker\Theta^p_*[g]$ if and only if $w\in\ker\Phi^p_*(g)$. First consider the case $g=e$. Then $w\in\ker\Phi^p_*(e)$ means
\[
\derev{}{t}{0}\Phi_{\exp(tw)}(p)=w_M(p)=0,
\]
and this is equivalent to $\Phi_{\exp(tw)}(p)=p$ for all $t\in\R$, since for any $s\in\R$ we have
\[
\derev{}{t}{s}\Phi_{\exp(tw)}(p)=\Phi_{\exp(sw)*}w_M(p).
\]
Hence $w\in\ker\Phi^p_*(e)$ if and only if $\exp(tw)\in G_p$ for all $t\in\R$, which implies $v=\pi_*w=0$, and this proves the assertion for $g=e$.

If $g\neq e$, then $w\in\ker\Phi^p_*(g)$ if and only if
\[
\Phi^p_*(w)= \derev{}{t}{0}\Phi^p(g\exp(tL_{g^{-1}*}w))=\Phi_{g*}(\Phi^p_*(L_{g^{-1}*}w))=0,
\]
which is equivalent to $L_{g^{-1}*}w\in \ker\Phi^p_*(e)$, since $\Phi_g$ is a diffeomorphism. Therefore, if $\Psi$ denotes the smooth action of $G$ on $G/G_p$ by left translations,
\[
v=\pi_*w=\derev{}{t}{0}\pi(g\exp(tL_{g^{-1}*}w))=\Psi_{g*}(\pi_*L_{g^{-1}*}w)=0. \qedhere
\]
\end{Proof}

From here we can deduce some useful characterizations of the isotropy Lie algebra $\gfrak_p:=T_eG_p$ of a point $p\in M$ and the tangent space of a $G$-orbit.

\begin{corollary} \label{Group-Actions:isotropy-algebra:p-tangent-orbit}
If $p\in M$ then
\[
\gfrak_p=\{\xi\in\gfrak:\xi_M(p)=0\} \quad \text{and} \quad T_p(G\cdot p)=\gfrak_M(p).
\]
\end{corollary}

\begin{Proof}
By the proof of \Cref{Group-Actions:injective-immersion}, if $\xi\in\gfrak$ is such that $\xi_M(p)=0$, then $\exp(t\xi)\in G_p$ for all $t\in\R$, and therefore $\xi\in\gfrak_p$. Conversely, if $\xi\in\gfrak_p$ then $\pi_*\xi=0$, so that
\[
\xi_M(p)=\Phi^p_*(\xi)=\Theta^p_*\circ\pi_*(\xi)=0.
\]
On the other hand, by the definition of infinitesimal generator, \cref{eq:def-infinitesimal-generator}, one has that $\Phi^p_*(\gfrak)=\gfrak_M(p)$. Since $\Phi_*^p=\Theta^p_*\circ \pi_*$ is surjective because both $\Theta^p_*$ and $\pi_*$ are surjective (when thinking of $\Theta^p$ as a diffeomorphism between $G/G_p$ and $G\cdot p$), the second claim follows.
\end{Proof}


\section{Moment Maps}

The concept of moment map is a generalization of the concept of momentum in physics, which appeared with the study of conserved quantities in mechanical systems. Let $(M,\omega)$ be a symplectic manifold and $G$ a Lie group acting symplectically on $M$ through $\Phi$ (that is, $\Phi_g$ is a symplectomorphism of $M$ on itself for each $g\in G$), and let $\gfrak$ be the Lie algebra of $G$.

\begin{definition}
We say that a smooth map $\mu:M\to \gfrak^*$ is a \textbf{moment map} for the action $\Phi$ if for every $\xi\in\gfrak$ we have that
\[
i_{\xi_M}\omega=d\muhat(\xi),
\]
where $\muhat:\gfrak\to \Cinf(M)$ is the comoment map, given by $\muhat(\xi)(p):=\mu(p)(\xi)$. In other words, every infinitesimal generator $\xi_M$ is a Hamiltonian vector field with energy $\muhat(\xi)$. If a moment map exists for the action $\Phi$ we say that the action is \textbf{Hamiltonian}.

We say that a moment map $\mu$ is $\Ad^*$-equivariant, or simply \textbf{equivariant}, if it is equivariant with respect to the coadjoint action on $\gfrak^*$, that is, if
\[
\mu\circ\Phi_g=\Ad^*(g)\circ\mu, \quad \text{for any $g\in G$.} \qedhere
\]
\end{definition}

Notice that $\mu:M\to\gfrak^*$ is smooth if and only if $\muhat(\xi)\in\Cinf(M)$ for every $\xi\in\gfrak$. Indeed, if we fix $\{\xi_i\}_i$ a basis for $\gfrak$ and $\{\lambda^i\}_i$ its dual basis for $\gfrak^*$, we can write $\mu=\mu_i\lambda^i$, for some $\mu_i:M\to\R$; then $\mu$ is smooth if and only if each $\mu_i$ is so, and it is immediate to see that $\mu_i=\muhat(\xi_i)$. Thus, the only condition for a moment map to exist is that every infinitesimal generator be a Hamiltonian vector field: if for every $\xi\in\gfrak$ there is some $f_\xi\in\Cinf(M)$ such that $\omega(\xi_M)=df_\xi$, since the Hamiltonian functions are only defined up to additive constants, the correspondence $\xi\mapsto f_\xi$ can be easily made linear on $\xi$, so that for each $p\in M$, the map $\mu(p):\gfrak\to\R$ defined by $\mu(p)(\xi):=f_\xi(p)$ lives in $\gfrak^*$, defining, hence, a moment map for the action.

As we said, the original interest in moment maps comes from their relationship with conserved quantities. This is expressed in the translation of Noether's classical theorem on symmetries and conserved quantities to the language of moment maps. We say that a Lie group $G$ acts by \textbf{symmetries} on a Hamiltonian system $(M,\omega,H)$ if it preserves $H$: if $\Phi_g^*H=H$ for every $g\in G$.

\begin{theorem}[Noether]\label{Moment-Maps:thm:Noether}
Let $(M,\omega,H)$ be a Hamiltonian system with a moment map $\mu$ arising from a Hamiltonian $G$-action $\Phi$ by symmetries. Then $\mu$ is a conserved quantity, that is, $\mu\circ\theta_t=\mu$, where $\theta_t$ is the flow of $X_H$.
\end{theorem}

\begin{Proof}
By the definition of comoment map, $\mu$ is a conserved quantity if and only if $\muhat(\xi)$ is a conserved quantity for every $\xi\in\gfrak$. On the one hand, we have that $\muhat(\xi)\circ\theta_0=\muhat(\xi)$, and on the other, that, at every point $p\in M$,
\begin{align*}
\der{}{t}(\muhat(\xi)\circ\theta_t)&=d\muhat(\xi)(X_H)=\omega(\xi_M,X_H)=-dH(\xi_M) \\
&= -\der{}{t}(H\circ\Phi_{\exp(t\xi)})=-\der{}{t}H=0.
\end{align*}
Hence $\muhat(\xi)\circ\theta_t=\muhat(\xi)$ and that ends the proof.
\end{Proof}

Before giving some concrete examples, we give a general strategy for constructing an equivariant moment map on exact symplectic manifolds. We say that $\theta$ is a potential for the symplectic structure $\omega$ if $\omega=-d\theta$.

\begin{proposition} \label{MomentMaps:construction}
Let $(M,\omega)$ be an exact symplectic manifold with potential $\theta$. Let $G$ be a Lie group acting smoothly on $M$ such that $\theta$ is invariant under the action: $\Phi_g^*\theta=\theta$ for every $g\in G$. Then the action is Hamiltonian with equivariant moment map given by the comoment map
\[
\muhat(\xi)=i_{\xi_M}\theta.
\]
\end{proposition} 

\begin{Proof}
The action is clearly symplectic because the differential commutes with any pullback. The map $\mu$ is obviously smooth since so are $\muhat(\xi)$ for any $\xi\in\gfrak$. Because $\theta$ is invariant under the $G$-action,
\[
\Lcal_{\xi_M}\theta=\derev{}{t}{0}\Phi_{\exp(t\xi)}^*\theta=0,
\]
and by Cartan's magic formula,
\[
d\muhat(\xi)=di_{\xi_M}\theta=\Lcal_{\xi_M}\theta-i_{\xi_M}d\theta=i_{\xi_M}\omega.
\]
Therefore, $\mu$ is a moment map for the action. To see that it is equivariant, we compute, using again the invariance of $\theta$ and also \cref{Group-Actions:inf-gen:properties:1} of \Cref{Group-Actions:inf-gen:properties}:
\begin{align*}
\Ad^*(g)\circ\mu(p)(\xi)&= \mu(p)(\Ad(g^{-1})\xi)=i_{(\Ad(g^{-1})\xi)_M}\theta(p) \\
&=i_{\Phi_{g^{-1}*}\xi_M}\theta(p)=\theta(\Phi_{g^{-1}*}(\xi_M(\Phi_g(p)))) \\
&=i_{\xi_M}\theta(\Phi_g(p))=\mu\circ\Phi_g(p). & \qedhere
\end{align*}
\end{Proof}

This construction can be particularized to the cotangent bundle.

\begin{corollary}\label{Moment-Maps:construction-corollary}
Let $M$ be a smooth manifold with a smooth $G$-action $\Phi$. Consider the $G$-action on $T^*M$ lifted from $\Phi$:
\[
\begin{array}{rrcl}
\Phi^\ell: & G\times T^*M & \longrightarrow & T^*M \\
& (g,\alpha) & \longmapsto & \Phi_{g^{-1}}^*(\alpha)
\end{array}.
\]
Then the lifted action acts by symplectomorphisms on the canonical symplectic structure on $T^*M$ and gives rise to an equivariant moment map defined by
\[
\mu(\alpha)(\xi)=\alpha(\xi_M(\pi(\alpha))),\quad \text{for $\alpha\in T^*M$,}
\]
where $\pi:T^*M\to M$ is the projection.
\end{corollary}

\begin{Proof}
First of all, for every $\alpha\in T_p^*M$ and $g\in G$,
\[
\pi\circ\Phi^\ell_g(\alpha)=\pi(\Phi^*_{g^{-1}}(\alpha))=\Phi_g(p)=\Phi_g\circ\pi(\alpha),
\]
which also implies that $\pi_*\xi_{T^*M}(\alpha)=\xi_M(\pi(\alpha))$. With this in mind, we can see that the tautological 1-form is invariant under the lifted action: for any $v\in T_\alpha T^*M$
\begin{align*}
(\Phi^{\ell*}_g\theta)_\alpha(v)&= \theta_{\Phi^*_{g^{-1}}(\alpha)}(\Phi^\ell_{g*}v)= \Phi^*_{g^{-1}}(\alpha)((\pi\circ\Phi^\ell_g)_*v) \\
&=\Phi^*_{g^{-1}}(\alpha)((\Phi_g\circ\pi)_*v)=\alpha(\pi_*v)=\theta_\alpha(v).
\end{align*}
This shows that the action indeed acts by symplectomorphisms and that, by \Cref{MomentMaps:construction}, the action is Hamiltonian with equivariant moment map
\[
\mu(\alpha)(\xi)=i_{\xi_{T^*M}}\theta(\alpha)=\theta_\alpha(\xi_{T^*M}(\alpha))=\alpha(\pi_*\xi_{T^*M}(\alpha))=\alpha(\xi_M(\pi(\alpha))). \qedhere
\]
\end{Proof}

The following properties of equivariant maps $\mu:M\to\gfrak^*$ will be very useful throughout the text.

\begin{proposition} \label{Moment-Maps:interesting-properties}
For any equivariant map $\mu:M\to \gfrak^*$ it is true that
\begin{enumerate}
	\item $\mu_*\xi_M(p)=\xi_{\gfrak^*}(\mu(p))$ for $p\in M$ and $\xi\in\gfrak$, so that $\mu_*(T_p(G\cdot p))=T_{\mu(p)}(G\cdot \mu(p))$, \label{moments:interesting-properties:1}
	\item $\xi_{\gfrak^*}(\alpha)=-\alpha\circ \ad\xi$ for $\alpha\in\gfrak^*$ and $\xi\in\gfrak$. \label{moments:interesting-properties:2}
\end{enumerate}
\end{proposition}

\begin{Proof}
\Cref{moments:interesting-properties:1} follows from the equivariance of $\mu$, and then, by \Cref{Group-Actions:isotropy-algebra:p-tangent-orbit},
\[
\mu_*(T_p(G\cdot p))=\{\mu_*\xi_M(p): \xi\in\gfrak\}=\{\xi_{\gfrak^*}(\mu(p)):\xi\in\gfrak\}=T_{\mu(p)}(G\cdot \mu(p)).
\]
\Cref{moments:interesting-properties:2} is a straightforward computation:
\[
\xi_{\gfrak^*}(\alpha)=\derev{}{t}{0}\Ad^*(\exp(t\xi))\alpha=\alpha \circ \derev{}{t}{0} \Ad(\exp(-t\xi)) = -\alpha\circ \ad\xi. \qedhere
\]
\end{Proof}

We now turn to the examples.

\begin{example}
Consider $(\R^3\times\R^3, \omega_0)$ acted upon by $\R^3$ by means of $a\cdot (q,p)=(q+a,p)$, where $a,q,p\in\R^3$. For a given $a\in\R^3$, the infinitesimal generator at $(q,p)\in\R^6$ is
\[
\derev{}{t}{0}(q+ta,p)=(a,p).
\]
Obviously the action preserves the symplectic potential $p_idq^i$. \Cref{MomentMaps:construction} then asserts that
\[
\mu(q,p)(a)=(i_{(a,p)}(p_idq^i))(q,p)=\prodesc{p}{a},
\]
with $\prodesc{\cdot}{\cdot}$ the usual inner product on $\R^3$, is an equivariant moment map for the action. Identifying $\R^3$ with its dual through this inner product, we obtain that the moment map is just $\mu(q,p)=p$. This is the linear momentum in classical mechanics.
\end{example}

\begin{example}\label{Moment-Maps:example:angular-momentum}
Consider $(\R^3\times\R^3, \omega_0)$ acted upon by $\SO(3)$ by means of $R\cdot (q,p)=(Rq,Rp)$, where $q,p\in\R^3$ and $R\in\SO(3)$. For a given $X\in\so(3)$, where we recall that $\so(3)$ is the set of traceless skew-symmetric real $3\times 3$ matrices, the infinitesimal generator at $(q,p)\in\R^6$ is
\[
\derev{}{t}{0}(e^{tX}q,e^{tX}p)=(Xq,Xp).
\]
This action also preserves the potential, since for any $R\in\SO(3)$,
\[
R^*(p_idq^i)_{(q,p)}(u,v)=(p_idq^i)_{(Rq,Rp)}(Ru,Rv)=\prodesc{Rp}{Ru}=\prodesc{p}{u}=(p_idq^i)_{(q,p)}(u,v),
\]
and the moment map is, then,
\[
\mu(q,p)(X)=(i_{(Xq,Xp)}(p_idq^i))(q,p)=\prodesc{p}{Xq}.
\]
The Lie algebra $\so(3)$ can be identified with $\R^3$ via
\begin{equation} \label{eq:identification-so(3)-R3}
\xi=(\xi_1,\xi_2,\xi_3)\longmapsto X_\xi:=\begin{pmatrix} 0 & -\xi_3 & \xi_2 \\ \xi_3 & 0 & -\xi_1 \\ -\xi_2 & \xi_1 & 0 \end{pmatrix},
\end{equation}
so that $X_\xi q=\xi\times q$, where $\times$ is the usual cross product on $\R^3$. With this identification, $\mu(q,p)(\xi)=\prodesc{p}{\xi\times q}=\prodesc{q\times p}{\xi}$, and with the further identification of $\R^3$ with its dual, we obtain $\mu(q,p)=q\times p$. This is the angular momentum in classical mechanics.
\end{example}

\begin{example}\label{Moment-Maps:example:left-trans}
Let $G$ be any Lie group acted upon by itself by left translations. For any $\xi\in\gfrak$, the infinitesimal generator at $g\in G$ is
\[
\xi_G(g)=\derev{}{t}{0}\exp(t\xi)g=R_{g*}\xi,
\]
where $R_g$ is right translation by $g$. By \Cref{Moment-Maps:construction-corollary}, the moment map of the lifted action to $T^*G$ is $\mu(\lambda)=\lambda\circ R_{g*}$ for any $\lambda\in T^*_gG$.

The cotangent bundle $T^*G$ can be trivialized through
\[
\begin{array}{rcl}
G\times\gfrak^* & \longrightarrow & T^*G \\
(g,\alpha) & \longmapsto & \alpha\circ R_{g^{-1}*}
\end{array},
\]
with inverse $\lambda\mapsto (g,\lambda\circ R_{g*})$ for $\lambda\in T_g^*G$. With this identification, the moment map is just projection onto the second factor:
\[
\mu(g,\alpha)=\mu(\alpha\circ R_{g^{-1}*})=\alpha. \qedhere
\]
\end{example}

We also present a couple of  examples on non-exact symplectic manifolds.

\begin{example}\label{ex:moment-map-s2}
Consider $(S^2,\omega)$, with $\omega$ the area form on $S^2$, acted upon by $\SO(3)$ by matrix multiplication. The form $\omega$ can be expressed in simple form as $\omega_x(u,v)=\prodesc{x}{u\times v}$, where $u,v\in T_xS^2$ are thought of as vectors in $\R^3$ orthogonal to $x$. For $X\in\so(3)$, the infinitesimal generator at $x\in S^2$ is
\[
\derev{}{t}{0}e^{tX}x=Xx.
\]
Then, using the identification of $\so(3)$ with $\R^3$ \eqref{eq:identification-so(3)-R3}, for any $u\in T_xS^2$ and $\xi\in\R^3$,
\[
(i_{X_\xi x}\omega)_x(u)=\prodesc{x}{X_\xi x\times u}=\prodesc{x}{(\xi\times x)\times u}=\prodesc{x\times (\xi\times x)}{u}.
\]
Using the identity $a\times (b\times c)=\prodesc{a}{c}b-\prodesc{a}{b}c$, we can rewrite this as
\[
(i_{X_\xi x}\omega)_x(u)=\prodesc{\prodesc{x}{x}\xi-\prodesc{x}{\xi}x}{u}=\prodesc{\xi}{u}=d(\prodesc{\xi}{\cdot})_x(u),
\]
so that $\mu(x)(\xi)=\prodesc{x}{\xi}$ defines a moment map for the action. Identifying again $\R^3$ with its dual we obtain that the moment map is just the inclusion $S^2\to \R^3$.
\end{example}

\begin{example}\label{ex:hamiltonian-as-moment-map}
Let $(M,\omega)$ be a symplectic manifold and $X_H$ a complete Hamiltonian vector field. The flow of $X_H$ gives rise to a Hamiltonian $\R$-action on $M$ with $H$ as moment map. Indeed, for any $a\in\R$ and $p\in M$, 
\[
a_M(p)=\derev{}{t}{0} \theta_{ta}(p)=aX_H(p), 
\]
so that $i_{a_M}\omega=ai_{X_H}\omega=adH=d\hat{H}(a)$.
\end{example}

Moment maps are not always submersions, as the next example shows.

\begin{example} \label{ex:moment-maps-are-not-submersions}
Let $\SO(3)$ act on $S^2\times S^2$ by the diagonal action: $R\cdot(x,y):=(Rx,Ry)$. By \Cref{ex:moment-map-s2}, the moment map for such an action (identifying $\so(3)^*$ with $\R^3$) is just $\mu(x,y)=x+y$, with image $\mu(S^2\times S^2)=\overline{B(0,2)}$, the closed ball of radius 2 and center 0. Hence, for any $x\in S^2$, since we can identify $T_xS^2$ with $\operatorname{span}(x)^\perp$, using the usual Euclidean metric on $\R^3$, we have that $\mu_*(x,-x)$ maps $T_xS^2\oplus T_{-x}S^2\cong \operatorname{span}(x)^\perp\oplus\operatorname{span}(x)^\perp$ to $\operatorname{span}(x)^\perp$, which is not all of $\R^3$.
\end{example}


\section{Poisson Geometry} 

The notion of a Poisson manifold is a generalization of the notion of a symplectic manifold. Its origins lie in the study of analytical mechanics, as did those of symplectic geometry. Poisson structures were introduced by A. Lichnerowicz in \cite{lichnerowicz} and have been extensively studied since then. 

\begin{definition}
A \textbf{Poisson manifold} is a pair $(M,\bracket{\cdot}{\cdot})$, where $M$ is a smooth manifold and $\bracket{\cdot}{\cdot}$ is a Lie bracket on the algebra $\Cinf(M)$ satisfying the Leibniz rule
\[
\bracket{f}{gh}=\bracket{f}{g}h+g\bracket{f}{h}, \quad\text{for all $f,g,h\in\Cinf(M)$.}
\]
If $(M',\bracket{\cdot}{\cdot}')$ is another Poisson manifold, we say that a diffeomorphism $F:M\to M'$ is \textbf{canonical} if the pullback $F^*$ preserves the brackets, i.e.,
\[
\bracket{F^*f}{F^*g}=F^*\bracket{f}{g}',\quad\text{for all $f,g\in\Cinf(M')$.}
\]
If such a canonical diffeomorphism exists we say that $M$ and $M'$ are canonically diffeomorphic. The bracket $\bracket{\cdot}{\cdot}$ is called \textbf{Poisson bracket}.
\end{definition}

The Leibniz rule asserts that $\bracket{f}{\cdot}$ is a derivation of $\Cinf(M)$, and so there is a vector field $X_f\in\Xfrak(M)$ such that $X_f=\bracket{f}{\cdot}$. Such a vector field is called the \textbf{Hamiltonian vector field} associated to $f$. As in the symplectic case, we say that $X\in\Xfrak(M)$ is Hamiltonian if $X=X_f$ for some $f\in\Cinf(M)$, and say that $f$ is a Hamiltonian or energy function for $X$.

Since $\bracket{\cdot}{\cdot}$ is bilinear, it is a derivation in both arguments, and so its value on $f,g\in\Cinf(M)$ depends only on the differentials $df$ and $dg$. That is, there is some bivector field $\Pi\in\Xfrak^2(M)$ such that
\begin{equation} \label{eq:bivector-field}
\bracket{f}{g}=\Pi(df,dg).
\end{equation}
Conversely, given a bivector field $\Pi\in\Xfrak^2(M)$, we can define a skew-symmetric bilinear bracket fulfilling the Leibniz rule by \cref{eq:bivector-field}. We say that a bivector field $\Pi\in\Xfrak^2(M)$ is Poisson if the associated bracket is Poisson, i.e., if it satisfies the Jacobi identity. It is apparent that it is enough to verify the Jacobi identity on functions whose differentials span $T_p^*M$ for each $p\in M$.

Poisson manifolds, hence, can be equivalently conceived as pairs $(M,\bracket{\cdot}{\cdot})$, where $\bracket{\cdot}{\cdot}$ is a Poisson bracket on $\Cinf(M)$, or as pairs $(M,\Pi)$, where $\Pi$ is a Poisson bivector field on $M$. Moreover, if $(M,\Pi)$ and $(M',\Pi')$ are Poisson manifolds, it is easy to see that a diffeomorphism $F:M\to M'$ is canonical if and only if $F_*\Pi=\Pi'$. There is a useful characterization of Poisson bivector fields in terms of the Schouten-Nijenhuis bracket. For a proof see \cite{vaisman}.

\begin{proposition}
A bivector field $\Pi\in\Xfrak^2(M)$ is Poisson if and only if $[\Pi,\Pi]=0$.
\end{proposition}

As we did in the symplectic case, consider the bundle morphism $\Pi:T^*M\to TM$ given by $\Pi(\lambda):=i_\lambda\Pi$. Since there are no nondegeneracy conditions on $\Pi$, this morphism may fail to be an isomorphism. The Hamiltonian vector field $X_f$ may hence be written as $X_f=\Pi(df)$.

We give now some basic examples of Poisson manifolds.

\begin{example} \label{PoissGeom:example:r2n}
Consider $\R^{2n}$ with Cartesian coordinates $(q^1,\dots,q^n,p_1,\dots,p_n,)$. Then the bracket defined by
\[
\bracket{f}{g}:=\derpar{f}{q^i}\derpar{g}{p_i}-\derpar{f}{p_i}\derpar{g}{q^i}
\]
is Poisson, and is called the standard Poisson structure on $\R^{2n}$. It is the Poisson bracket one usually encounters in classical mechanics on $\R^{2n}$ (see {\cite[\S 42]{landau-lifshitz}}). The corresponding Poisson bivector field is obviously
\[
\Pi=\derpar{}{q^i}\wedge\derpar{}{p_i}. \qedhere
\]
\end{example}

\begin{example} \label{ex:symplectic-are-poisson}
Let $M$ be a smooth manifold and $\omega\in\Omega^2(M)$ nondegenerate. We define the bracket of $f,g\in\Cinf(M)$ as
\[
\bracket{f}{g}:=-\omega(X_f,X_g)=dg(X_f)=X_fg.
\]
We claim that it is Poisson if and only if $d\omega=0$. In particular if $(M,\omega)$ is a symplectic manifold, then it is also a nondegenerate Poisson manifold (meaning that its Poisson bivector field is nondegenerate).

The bracket is evidently bilinear and alternating. If $f,g,h\in\Cinf(M)$, by the Leibniz rule for differentials,
\[
\bracket{f}{gh}= d(gh)(X_f)=hdg(X_f)+gdh(X_f)=h\bracket{f}{g}+g\bracket{f}{h}.
\]
Lastly, by the formula for the exterior differential of a 2-form, we have that
\begin{align*}
d\omega(X_f,X_g,X_h)&=\begin{multlined}[t] X_f\omega(X_g,X_h)+X_g\omega(X_h,X_f)+X_h\omega(X_f,X_g)\\
-\omega([X_f,X_g],X_h)-\omega([X_g,X_h],X_f)-\omega([X_h,X_f],X_g) \end{multlined}\\
&=\begin{multlined}[t]-\bracket{f}{\bracket{g}{h}}-\bracket{g}{\bracket{h}{f}}-\bracket{h}{\bracket{f}{g}} \\
-\omega([X_f,X_g],X_h)-\omega([X_g,X_h],X_f)-\omega([X_h,X_f],X_g). \end{multlined}
\end{align*}
By Cartan's magic formulas, and taking into account that $di_{X_f}\omega=d^2f=0$ (and so $\Lcal_{X_f}\omega=i_{X_f}d\omega$),
\begin{equation}\label{eq:differential-of-bracket-symp}
\begin{aligned}
d\bracket{f}{g}&= d(\omega(X_g,X_f))=d(i_{X_f}i_{X_g}\omega)=\Lcal_{X_f}(i_{X_g}\omega)-i_{X_f}(di_{X_g}\omega) \\
&=\Lcal_{X_f}(i_{X_g}\omega)=i_{[X_f,X_g]}\omega+i_{X_g}\Lcal_{X_f}\omega\\
&=i_{[X_f,X_g]}\omega+i_{X_g}i_{X_f}d\omega.
\end{aligned}
\end{equation}
Hence,
\begin{align*}
\omega([X_f,X_g],X_h)&= i_{[X_f,X_g]}\omega(X_h)=d\bracket{f}{g}(X_h)-i_{X_g}i_{X_f}d\omega(X_h) \\
&=\bracket{h}{\bracket{f}{g}}-d\omega(X_f,X_g,X_h),
\end{align*}
and similarly for $\omega([X_g,X_h],X_f)$ and $\omega([X_h,X_f],X_g)$. Therefore
\[
d\omega(X_f,X_g,X_h)=-2\bracket{f}{\bracket{g}{h}}-2\bracket{g}{\bracket{h}{f}}-2\bracket{h}{\bracket{f}{g}}+3d\omega(X_f,X_g,X_h),
\]
and so,
\[
d\omega(X_f,X_g,X_h)=\bracket{f}{\bracket{g}{h}}+\bracket{g}{\bracket{h}{f}}+\bracket{h}{\bracket{f}{g}}.
\]
If $(U,(x^i))$ is a local chart on $M$, then $X\in\Xfrak(U)$ can be uniquely written as $X=X^iX_{x^i}$ for some functions $X^i\in\Cinf(U)$ (indeed, if we write $\omega^{ij}$ for the components of the inverse matrix of $(\omega_{ij})$, then $X_{x^i}=\omega^{ij}\partial_j$). Thus, it is clear that $\bracket{\cdot}{\cdot}$ fulfills the Jacobi identity if and only if $d\omega=0$.

We conclude that a symplectic manifold can be equivalently considered as a pair $(M,\omega)$, with $\omega$ a nondegenerate closed 2-form, or as a pair $(M,\Pi)$, with $\Pi$ a nondegenerate Poisson bivector field. The corresponding bundle isomorphisms are related by $\Pi=\omega^{-1}$. In addition, a diffeomorphism $F\in\Diff(M)$ is canonical if and only if it is a symplectomorphism, since if $F$ is a symplectomorphism or canonical, it is easily seen that $X_{F^*f}=F^{-1}_*X_f$, which immediately implies the other.
\end{example}

If $(M,\omega)$ is a symplectic manifold carrying a Hamiltonian action of some Lie group $G$ with equivariant moment map $\mu$, then, because of \cref{Group-Actions:inf-gen:properties:2} of \Cref{Group-Actions:inf-gen:properties}, \cref{eq:differential-of-bracket-symp} implies that $d\{\muhat(\xi),\muhat(\eta)\}=-d\muhat([\xi,\eta])$, i.e., $\{\muhat(\xi),\muhat(\eta)\}+\muhat([\xi,\eta])$ is a constant function. By \Cref{Moment-Maps:interesting-properties}, this constant is, for any $p\in M$,
\begin{align*}
\{\muhat(\xi),\muhat(\eta)\}(p)+\muhat([\xi,\eta])(p)& = -\omega_p(\xi_M(p),\eta_M(p))+\mu(p)\circ\ad\xi(\eta) \\
&=d\muhat(\eta)(\xi_M(p))+\mu(p)\circ\ad\xi(\eta) \\
&=(\mu_*\xi_M(p)+\mu(p)\circ\ad\xi)(\eta)=0.
\end{align*}

Hence, in fact we have that $\{\muhat(\xi),\muhat(\eta)\}=-\muhat([\xi,\eta])$, that is, the comoment map is a Lie algebra antihomomorphism between $(\gfrak,[\cdot,\cdot])$ and $(\Cinf(M),\bracket{\cdot}{\cdot})$. Conversely, it can be seen that if the comoment map is a Lie algebra antihomomorphism, then the moment map is equivariant with respect to the elements of the identity component of $G$ (see {\cite[Lem. 5.16]{mcduff-salamon}}). It can also be seen that the obstruction for a moment map to be equivariant is of cohomological type, concerning the vanishing of a particular cocycle class in the Lie algebra cohomology of $\gfrak$ (see {\cite[Lem 5.15]{mcduff-salamon}} and comments afterwards).

Let $G$ be a Lie group acting smoothly through $\Phi$ on a Poisson manifold $(M,\bracket{\cdot}{\cdot})$. We say that the action is \textbf{canonical} if $\Phi_g$ is a canonical diffeomorphism for every $g\in G$. Furthermore, as in the symplectic case, we say that the action is Hamiltonian if there is a moment map for the action, meaning a smooth map $\mu:M\to\gfrak^*$ such that the comoment map satisfies
\[
\Pi(d\muhat(\xi))=\xi_M,\quad\text{ for every $\xi\in\gfrak$.}
\]
The underlying idea is the same as in the symplectic case: for each $\xi\in\gfrak$, the vector field $\xi_M$ is Hamiltonian with $\muhat(\xi)$ as energy function.

For further details on Poisson geometry consult \cite{vaisman, weinstein-poisson}.

\chapter{Reduction}\label{chap:reduction}

It is a common practice in physics to ``divide out symmetries'' when trying to solve the dynamics of a mechanical system. For instance, if $(\R^3\times\R^3,\omega_0,H)$ is a Hamiltonian system such that $H$ only depends on $\norm{q}$ and $\norm{p}$, where $(q,p)\in \R^3\times\R^3$, then it is invariant under the action of $\SO(3)$ described in \Cref{Moment-Maps:example:angular-momentum} and the angular momentum $L(q,p):=q\times p$ gives rise to three conserved quantities by Noether's theorem. Hence, we can limit our study of the dynamics to $L^{-1}(\xi)$, for any $\xi\in\R^3$ such that $L^{-1}(\xi)$ defines a smooth manifold. Moreover, there is still some symmetry left on the level set $L^{-1}(\xi)$, meaning that there is a subgroup of $\SO(3)$ still acting on $L^{-1}(\xi)$. Therefore, we can quotient $L^{-1}(\xi)$ by this subgroup, identifying all the points in the same orbit, and obtain a ``reduced'' space, which will have fewer degrees of freedom than the original one, and which will presumably be easier to solve. We can then ``lift'' the solved motion to the unreduced space and obtain a motion therein.

In this chapter we study, under some regularity conditions, the reduction of manifolds by symmetry groups in both the symplectic and the Poisson cases.


\section{Symplectic Reduction}\label{sect:symplectic-reduction}

\subsection{Marsden-Weinstein Reduction}

Symplectic reduction was introduced by Marsden and Weinstein in 1974 in a short paper \cite{marsden-weinstein}. 

First of all, we recall that a $G$-action on a manifold $M$, where $G$ is some Lie group, is said to be \textbf{free} if the isotropy group of any element is trivial. Hence, by \Cref{Group-Actions:injective-immersion}, each orbit is an immersed submanifold diffeomorphic to $G$. The action is said to be \textbf{proper} if the map $(g,p)\mapsto (g,\Phi_g(p))$, for $g\in G$ and $p\in M$, is proper, in the sense that the preimage of compact sets are compact sets. It is well known (see for instance {\cite[Thm. 21.10]{lee}}) that if the action is both free and proper, then $M/G$ has a unique smooth structure such that the projection $\pi: M\to M/G$ is a submersion. In this case $T_{\pi(p)}(M/G)\cong T_pM/T_p(G\cdot p)$.

Suppose now that $(M,\omega)$ is a symplectic manifold, and that the action is Hamiltonian with associated equivariant moment map $\mu:M\to\gfrak^*$. We say that $\alpha\in\gfrak^*$ is a \textbf{clean value} of $\mu$ if $\mu^{-1}(\alpha)$ is a submanifold such that $T_p\mu^{-1}(\alpha)=\ker \mu_*(p)$. For instance, if $\alpha$ is a regular value, then it is evidently a clean value.

We first prove a useful technical lemma. By $^\circ$ we denote the annihilator of a vector subspace.

\begin{lemma} \label{lemma:kermu:symplectic}
\begin{enumerate}
	\item For any $p\in M$, $\ker \mu_*(p)=\gfrak_M(p)^\perp$ and $\im \mu_*(p)=\gfrak_p^\circ$. \label{lemma:kermu:symplectic:1}
	\item If $\alpha\in\gfrak^*$ is a clean value for $\mu$ and $p\in\mu^{-1}(\alpha)$, then $T_p(G_\alpha\cdot p)=T_p(G\cdot p)\cap T_p\mu^{-1}(\alpha)$. \label{lemma:kermu:symplectic:2}
\end{enumerate}
\end{lemma}

\begin{Proof}
\Cref{lemma:kermu:symplectic:1} follows from
\[
\omega(\xi_M(p),v)=d\muhat(\xi)(v)=\mu_*v(\xi), \quad \text{ for $v\in T_pM$ and $\xi\in\gfrak$,}
\]
and the nondegeneracy of $\omega$.

To see \cref{lemma:kermu:symplectic:2}, since $\alpha$ is a clean value, $\xi_M(p)\in T_p\mu^{-1}(\alpha)$ if and only if $\mu_*\xi_M(p)=0$. By \Cref{Moment-Maps:interesting-properties}, this is precisely $\xi_{\gfrak^*}(\alpha)=0$, which in turn is equivalent, by \Cref{Group-Actions:isotropy-algebra:p-tangent-orbit}, to $\xi\in\gfrak_\alpha$. Hence $\xi_M(p)\in(\gfrak_\alpha)_M(p)=T_p(G_\alpha\cdot p)$.  
\end{Proof}

To motivate the construction of the symplectic reduction, consider $\alpha$ a clean value of $\mu$. For any $p\in\mu^{-1}(\alpha)$, by \cref{lemma:kermu:symplectic:1} of \Cref{lemma:kermu:symplectic}, $T_p\mu^{-1}(\alpha)=\gfrak_M(p)^\perp$. This means that, if $i:\mu^{-1}(\alpha)\to M$ is the inclusion, the elements of $\gfrak_M(p)\cap T_p\mu^{-1}(\alpha)=T_p(G_\alpha\cdot p)$ lie in the kernel of the 2-form $i^*\omega$. One would think, then, that by dividing by $G_\alpha$ we will obtain a symplectic form on $\mu^{-1}(\alpha)$, and this is precisely what happens. This makes sense since $\mu^{-1}(\alpha)$ is $G_\alpha$-invariant by the equivariance of $\mu$.

\begin{theorem}[Marsden-Weinstein]\label{Reduction:thm:Marsden-Weinstein}
Let $(M,\omega)$ be a symplectic manifold with an equivariant moment map $\mu$ arising from a Hamiltonian $G$-action $\Phi$. Let $\alpha\in\gfrak^*$ be a clean value of $\mu$ and let $i:\mu^{-1}(\alpha)\to M$ be the inclusion. Suppose that $G_\alpha$ acts freely and properly on $\mu^{-1}(\alpha)$. Then there is a unique symplectic structure $\omega_\alpha$ on $M\dbar_\alpha G:=\mu^{-1}(\alpha)/G_\alpha$ such that $i^*\omega=\pi^*\omega_\alpha$, where $\pi: \mu^{-1}(\alpha)\to M\dbar_\alpha G$ is the projection.
\end{theorem}

\begin{Proof}
If $\omega_\alpha$ exists, then its value is wholly determined by the relation $i^*\omega=\pi^*\omega_\alpha$, so uniqueness is clear. To see existence, for any $\pi_*u,\pi_*v\in T_{\pi(p)}(M\dbar_\alpha G)$, with $u,v\in T_p\mu^{-1}(\alpha)$, define
\[
\omega_\alpha(\pi_*u,\pi_*v):=i^*\omega(u,v).
\]
To see that it is well defined, assume that $u'\in T_p\mu^{-1}(\alpha)$ is such that $\pi_*u'=\pi_*u$. Then $u-u'\in\ker\pi_*=T_p(G_\alpha\cdot p)$, and by the comments preceding this theorem, $u-u'\in\ker i^*\omega$ and, hence,
\[
i^*\omega(u,v)-i^*\omega(u',v)=0, \text{ for all $v\in T_p\mu^{-1}(\alpha)$.}
\]
On the other hand, suppose that $q\in \mu^{-1}(\alpha)$ is such that $\pi(q)=\pi(p)$, i.e., there is some $g\in G_\alpha$ such that $p=\Phi_g(q)$, and that there are $u',v'\in T_q\mu^{-1}(\alpha)$ with $\pi_*u'=\pi_*u$ and $\pi_*v'=\pi_*v$. Then $\Phi_{g*}u',\Phi_{g*}v'\in T_p\mu^{-1}(\alpha)$ are such that $\pi_*\Phi_{g*}u'=\pi_*u$ and $\pi_*\Phi_{g*}v'=\pi_*v$ because $\pi\circ\Phi_g=\pi$, and therefore
\[
i^*\omega(\Phi_{g*}u',\Phi_{g*}v')=i^*\omega(u,v).
\]
Since the action is symplectic and $\Phi_g$ commutes with $i$, we finally obtain that
\[
i^*\omega(u',v')=i^*\omega(u,v),
\]
and $\omega_\alpha$ is well defined.

By the previous comments, it is clear that $\omega_\alpha$ is nondegenerate. Because the differential commutes with pullbacks, we have that $\pi^*d\omega_\alpha=i^*d\omega=0$. Because $\pi$ is a submersion, this is equivalent to $d\omega_\alpha=0$, and so $\omega_\alpha$ indeed defines a symplectic structure on $M\dbar_\alpha G$.
\end{Proof}

It is useful to have the following diagram in mind
\[
	\begin{tikzcd}
  	(\mu^{-1}(\alpha),i^*\omega=\pi^*\omega_\alpha) \arrow[r,"i"] \arrow[d,"\pi",swap] & (M,\omega) \\
	(M\dbar_\alpha G,\omega_\alpha)
	\end{tikzcd}
\]

The manifold $M\dbar_\alpha G$ is called the \textbf{symplectic reduced space}. Of course, everything works fine because we are in the most friendly scenario for symplectic reduction. As shown in \cite{lerman-tolman}, if the action is not taken to be free, then $M\dbar_\alpha G$ is not a manifold, but an orbifold (see this same paper for the concept of orbifold, or look at \cite{satake-2}). If $\alpha$ is a singular value, then $M\dbar_\alpha G$ was shown in \cite{sjamaar-lerman} to be a symplectic stratified space, which, roughly speaking, is a collection of symplectic manifolds that fit together nicely (the theorem stating that the reduced space decomposes as a disjoint union of symplectic manifolds is presented as \Cref{thm:sjamaar-lerman} in \Cref{sect:cross-section-theorems}). The book \cite{ortega-ratiu} is a comprehensive exposition of reduction theory in general symplectic and Poisson manifolds, with no regularity conditions whatsoever. For our purposes, though, the regular Marsden-Weinstein theorem is enough.

One of the appeals of symplectic reduction is that by taking a simple symplectic manifold $M$ and reducing it with respect to some action we can obtain symplectic structures on relatively complicated objects for which a direct presentation of its symplectic structure would seem rather magical or special. We illustrate this fact with two examples.

\begin{example}
Here we follow {\cite[Ex. 4.3.4(iv)]{abraham-marsden}}. Let $X_H$ be a complete Hamiltonian vector field on a symplectic manifold $(M,\omega)$. By \Cref{ex:hamiltonian-as-moment-map}, the flow of $X_H$ gives rise to a Hamiltonian $\R$-action on $M$ with $H$ as moment map. Hence, if $a\in\R$ is a clean value of $H$ we obtain that $H^{-1}(a)/\R$ is a symplectic manifold if the action is free and proper on $H^{-1}(a)$.

Take $(M,\omega)=(\R^{2n},\omega_0)$, the standard symplectic structure, and take the harmonic oscillator Hamiltonian
\[
H(q,p)=\frac{1}{2}\sum_{i=1}^n ((q^i)^2+p_i^2).
\]
Then by \Cref{SympGeom:example:r2n} we have that
\[
X_H=p_i\derpar{}{q^i}-q^i\derpar{}{p_i},
\]
whose flow is
\[
\theta_t(q,p)=(q\cos t+p\sin t, p\cos t-q\sin t).
\]
Since $\theta_t$ is $2\pi$-periodic in $t$, it defines a symplectic action of $S^1$ on $M$, which is obviously free. Since $S^1$ is compact, the action is proper, and since $1/2$ is a regular value of $H$, 
\[
H^{-1}(1/2)/\R=H^{-1}(1/2)/S^1=S^{2n-1}/S^1=\C P^{n-1}
\]
is a $2(n-1)$-dimensional symplectic manifold.
\end{example}

Before giving the second example, we need the expression for the canonical symplectic form on the cotangent bundle of a Lie group.

\begin{proposition} \label{prop:symplectic-form-cotangent-bundle-lie-group}
Let $G$ be a Lie group and $\gfrak$ its Lie algebra. Identify $T^*G$ with $G\times\gfrak^*$ via right translations and let $\omega$ be the canonical symplectic form on $T^*G$. Then, for $(g,\alpha)\in G\times\gfrak^*$, $\xi,\eta\in \gfrak$ and $\beta,\gamma\in\gfrak^*$,
\[
\omega(R_{g*}\xi+\beta,R_{g*}\eta+\gamma)=\gamma(\xi)-\beta(\eta)-\alpha([\xi,\eta]).
\]
\end{proposition}

\begin{Proof}
The tautological form $\theta\in\Omega^1(T^*G)$ is given by $\theta_{(g,\alpha)}(R_{g*}\xi+\beta)=(g,\alpha)(R_{g*}\xi)=\alpha(\xi)$. Consider the $G$-action on itself by left translations. We already saw in \Cref{Moment-Maps:example:left-trans} that for any $\zeta\in\gfrak$, we have that $\zeta_G(g)=R_{g*}\zeta$. Hence, the vector field $X_{\xi,\beta}:=\xi_G+\beta$ is such that at $(g,\alpha)$ its value is $R_{g*}\xi+\beta$, and similarly for $X_{\eta,\gamma}:=\eta_G+\gamma$ and $R_{g*}\eta+\gamma$. Note that the flow of $X_{\xi,\beta}$ is
\begin{equation} \label{eq:flow-Xxibeta}
\phi_t^{\xi,\beta}(h,\lambda)=(\exp(t\xi)h,\lambda+t\beta), \quad\text{ for $(h,\lambda)\in G\times\gfrak^*$.}
\end{equation}
By the formula for the differential of a 1-form, we have that
\begin{align*}
\omega(R_{g*}\xi+\beta,R_{g*}\eta+\gamma)&=-d\theta(R_{g*}\xi+\beta,R_{g*}\eta+\gamma)=-d\theta(X_{\xi,\beta},X_{\eta,\gamma})(g,\alpha) \\
&= \big(-X_{\xi,\beta}(\theta(X_{\eta,\gamma}))+X_{\eta,\gamma}(\theta(X_{\xi,\beta}))+\theta([X_{\xi,\beta},X_{\eta,\gamma}])\big) (g,\alpha).
\end{align*}
It is not hard to see that $\theta(X_{\eta,\gamma})(h,\lambda)=\lambda(\eta)$, for $(h,\lambda)\in G\times\gfrak^*$, so that $X_{\xi,\beta}(\theta(X_{\eta,\gamma}))(g,\alpha)=\beta(\eta)$. Similarly, $X_{\eta,\gamma}(\theta(X_{\xi,\beta}))(g,\alpha)=\gamma(\xi).$ Using \cref{eq:flow-Xxibeta} and that the Lie bracket can be computed as the Lie derivative, one can also see that 
\[
[X_{\xi,\beta},X_{\eta,\gamma}](h,\lambda)=-R_{h*}([\xi,\eta]),\quad\text{ for $(h,\lambda)\in G\times\gfrak^*$,}
\]
and therefore $\theta([X_{\xi,\beta},X_{\eta,\gamma}])(g,\alpha)=-\alpha([\xi,\eta])$. Hence the result.
\end{Proof}

\begin{example}
We retake \Cref{Moment-Maps:example:left-trans}. The moment map was given by $\mu(\lambda)=\lambda\circ R_{g*}$ for $\lambda\in T_g^*G$. Through the identification of $T^*G$ with $G\times\gfrak^*$, the moment map could be written as $\mu(g,\alpha)=\alpha$, which shows that any value $\alpha\in\gfrak^*$ is regular. The lifted action to the cotangent bundle can be written as
\[
h\cdot(g,\alpha)=(hg,\Ad^*(h)\alpha).
\]

Also, $\mu^{-1}(\alpha)=G\times\{\alpha\}$, so that $T^*G\dbar_\alpha G\cong G/G_\alpha$, where $G/G_\alpha$ here is the orbit space by left translations, which equals the set of right cosets. Right cosets are in bijective correspondence with left cosets by the map $G_\alpha g\mapsto g^{-1}G_\alpha$, so that, by \Cref{Group-Actions:injective-immersion}, in fact $T^*G\dbar_\alpha G\cong G\cdot\alpha$. To see that the action is proper, by {\cite[Prop. 21.5]{lee}}, it suffices to see that if $(g_i)$ and $(h_i)$ are sequences in $G$ such that $(g_i)$ and $(g_ih_i)$ converge, then a subsequence of $(h_i)$ converges. But this is obvious, since if $g:=\lim_ig_i$ and $k:=\lim_i g_ih_i$, then, by continuity of the product and inverse maps on $G$, $\lim_ ih_i=\lim_i g_i^{-1}(g_ih_i)=g^{-1}k$. On the other hand, the action is obviously free. By \Cref{Reduction:thm:Marsden-Weinstein}, $G\cdot\alpha$ has a symplectic structure.

Let us compute this structure explicitly. Because of the identification of right cosets with left cosets, the projection $\pi:\mu^{-1}(\alpha)\to G\cdot\alpha$ is given by $\pi(g,\alpha)=g^{-1}\cdot\alpha=\Ad^*(g^{-1})\alpha$. Hence, for any $\xi\in\gfrak$, viewing $R_{g*}\xi\in T_gG$ as an element of $T_{(g,\alpha)}\mu^{-1}(\alpha)$,
\begin{align*}
\pi_*R_{g*}\xi&=\derev{}{t}{0}\pi(\exp(t\xi)g,\alpha)=\derev{}{t}{0}\Ad^*(g^{-1}\exp(-t\xi)g)\Ad^*(g^{-1})\alpha \\
&=-\Ad^*(g^{-1})\alpha\circ\ad(C_{g^{-1}*}\xi) =(\Ad(g^{-1})\xi)_{\gfrak^*}(g^{-1}\cdot\alpha),
\end{align*}
where the last equality follows from \Cref{Moment-Maps:interesting-properties}.

By \Cref{Reduction:thm:Marsden-Weinstein}, then, 
\[
\omega_\alpha((\Ad(g^{-1})\xi)_{\gfrak^*}(g^{-1}\cdot\alpha),(\Ad(g^{-1})\eta)_{\gfrak^*}(g^{-1}\cdot\alpha))=\omega(R_{g*}\xi,R_{g*}\eta),\text{ for any $\xi,\eta\in T_gG$,}
\]
where $\omega$ is the canonical symplectic form on $T^*G$. By \Cref{prop:symplectic-form-cotangent-bundle-lie-group}, $\omega(R_{g*}\xi,R_{g*}\eta)=-\alpha([\xi,\eta])$. 

We conclude that for any $\xi,\eta\in\gfrak$ and $\beta=g\cdot\alpha$ for $g\in G$, since $[\Ad(g)\xi,\Ad(g)\eta]=\Ad(g)[\xi,\eta]$,
\[
\omega_\alpha(\xi_{\gfrak^*}(\beta), \eta_{\gfrak^*}(\beta))=-\beta([\xi,\eta]).
\]

This symplectic structure on the coadjoint orbit $G\cdot\alpha$ is commonly known as the Kirillov-Kostant-Souriau symplectic form, first proven to be a symplectic form by Kostant {\cite[Thm. 5.3.1]{kostant}}. By means of Marsden-Weinstein reduction this structure arises naturally. This example can be found in \cite{marsden-weinstein}.
\end{example}

\subsection{Reduced Symplectic Dynamics}

As has already been said, the interest in reduction is to simplify a given system. So far we have reduced the symplectic manifold---the phase space---, but from the dynamical point of view we are also interested in how the dynamics on the manifold behave under the reduction process. The answer is that we obtain reduced dynamics in the reduced system, related in a natural way to the dynamics on the original system.

\begin{theorem}[Reduced symplectic dynamics] \label{thm:reduced-symplectic-dynamics}
Let $(M,\omega,H)$ be a Hamiltonian system with moment map $\mu$ arising from a Hamiltonian $G$-action $\Phi$ by symmetries. Let $\alpha\in\gfrak^*$ be a clean value of $\mu$ such that $G_\alpha$ acts freely and properly on $\mu^{-1}(\alpha)$. Then the flow of $X_H$ leaves $\mu^{-1}(\alpha)$ invariant and induces a Hamiltonian flow on $M\dbar_\alpha G$ with Hamiltonian $h$ given by $i^*H=\pi^*h$, where $i:\mu^{-1}(\alpha)\to M$ is the inclusion and $\pi:\mu^{-1}(\alpha)\to M\dbar_\alpha G$ the projection. In particular, $X_H$ and $X_h$ are $\pi$-related. The function $h$ is called \textbf{reduced Hamiltonian} or energy.
\end{theorem}

\begin{Proof}
By Noether's theorem (\Cref{Moment-Maps:thm:Noether}), $\mu$ is a conserved quantity, so that the flow of $X_H$, say $\Theta_t$, preserves $\mu^{-1}(\alpha)$. Hence, it gives rise to a well-defined flow $\theta_t$ on $M\dbar_\alpha G$ such that $\pi\circ\Theta_t=\theta_t\circ\pi$. Let $X$ be the vector field associated to $\theta_t$. Then $X_H$ and $X$ are $\pi$-related: for any $p\in\mu^{-1}(\alpha)$,
\[
X\circ\pi(p)=\derev{}{t}{0}\theta_t\circ\pi(p)=\derev{}{t}{0}\pi\circ\Theta_t(p)=\pi_*X_H(p).
\]
The invariance of $H$ under the $G$-action ensures the existence of a smooth function $h$ in $M\dbar_\alpha G$ uniquely determined by $i^*H=\pi^*h$. For any $p\in\mu^{-1}(\alpha)$ and $v\in T_p\mu^{-1}(\alpha)$, we have that
\begin{align*}
(i_X\omega_\alpha)_{\pi(p)}(\pi_*v)&= \omega_\alpha(X\circ\pi(p),\pi_*v)=\omega_\alpha(\pi_*X_H(p),\pi_*v)\\
&=i^*\omega(X_H(p),v)=i^*dH(v)=\pi^*dh(v)=dh(\pi_*v).
\end{align*}
Thus, $X=X_h$, and we are done.
\end{Proof}

The ultimate interest in the dynamics is, of course, to recover the motion of the original system (recall that a motion of the Hamiltonian system $(M,\omega,H)$ is an integral curve of $X_H$) from the motion of the reduced system. We follow \cite{abraham-marsden} in this exposition.

\begin{lemma} \label{Reduction:lemma:diff-of-action}
Let $M$ be a manifold with a smooth $G$-action $\Phi$, and let $p\in M$ and $g\in G$. For $v\in T_pM$ and $\xi\in\gfrak$, and identifying $T_{(g,p)}(G\times M)\cong T_gG\oplus T_pM$ in the natural way, we have that
\[
\Phi_*(R_{g*}\xi+v)=\Phi_{g*}\left((\Ad(g^{-1})\xi)_M(p)+v\right).
\]
\end{lemma}

\begin{Proof} By the Leibniz rule for pushforwards,
\[
\Phi_*(R_{g*}\xi+v)= \Phi^p_*R_{g*}\xi+\Phi_{g*}v=\derev{}{t}{0}\Phi^p(g\exp(tC_{g^{-1}*}\xi))+\Phi_{g*}v=\Phi_{g*}\left((\Ad(g^{-1})\xi)_M(p)+v\right). \qedhere
\]
\end{Proof}

\begin{theorem}[Lifted motion]
Let $(M,\omega,H)$ be a Hamiltonian system with moment map $\mu$ arising from a Hamiltonian $G$-action $\Phi$ by symmetries. Let $\alpha\in\gfrak^*$ be a clean value of $\mu$ such that $G_\alpha$ acts freely and properly on $\mu^{-1}(\alpha)$, and let $(M\dbar_\alpha G, \omega_\alpha,h)$ be the reduced system. Let $\gamma:I\to M\dbar_\alpha G$ be a motion of the reduced system and $\beta:I\to M$ any smooth curve such that $\pi\circ\beta=\gamma$. Let $\xi:I\to\gfrak_\alpha$ be a smooth curve such that
\begin{equation} \label{eq:recover-dynamics-symplectic:algebraic}
\xi(t)_M(\beta(t))=X_H(\beta(t))-\betadot(t)
\end{equation}
and $g:I\to G$ a smooth curve such that
\begin{equation} \label{eq:recover-dynamics-symplectic:differential}
\gdot(t)=L_{g(t)*}\xi(t).
\end{equation}
Then the curve defined by $\Gamma(t):=\Phi_{g(t)}(\beta(t))$ is a motion of the original system. The fact that $\gamma$ is a motion for the reduced system ensures the existence and uniqueness of the solution to \cref{eq:recover-dynamics-symplectic:algebraic}.
\end{theorem}

\begin{Proof}
We compute, using \Cref{Reduction:lemma:diff-of-action},
\begin{align*}
\Gammadot(t)&=\Phi_*(\gdot(t)+\betadot(t)) \\
&=\Phi_{g(t)*}(\betadot(t)+(L_{g(t)^{-1}*}\gdot(t))_M(\beta(t))) \\
&=\Phi_{g(t)*}X_H(\beta(t)).
\end{align*}
Since $G$ acts symplectically by symmetries, for any $p\in M$ and $v\in T_{\Phi_g(p)}M$, 
\begin{align*}
\omega_{\Phi_g(p)}(X_H(\Phi_g(p)),v)&=dH(v)=d(H\circ\Phi_{g^{-1}})(v)=dH(\Phi_{g^{-1}*}v) \\
&=\omega_p(X_H(p),\Phi_{g^{-1}*}v)=(\Phi^*_g\omega)_p(X_H(p),\Phi_{g^{-1}*}v)\\
&=\omega_{\Phi_g(p)}(\Phi_{g*}X_H(p),v),
\end{align*}
so that $\Phi_{g*}X_H(p)=X_H\circ\Phi_g(p)$. Hence,
\[
\Gammadot(t)=\Phi_{g(t)*}X_H(\beta(t))=X_H(\Phi_{g(t)}(\beta(t)))=X_H(\Gamma(t)).
\]

To see the existence and uniqueness of the solution to \cref{eq:recover-dynamics-symplectic:algebraic}, note that, since $X_H$ and $X_h$ are $\pi$-related by \Cref{thm:reduced-symplectic-dynamics}, we have that
\[
\pi_*(X_H(\beta(t))-\betadot(t))=X_h(\gamma(t))-\gammadot(t)=0.
\]
Therefore $X_H(\beta(t))-\betadot(t)\in\ker\pi_*(\beta(t))=T_{\beta(t)}(G_\alpha\cdot\beta(t))=(\gfrak_\alpha)_M(\beta(t))$ for each $t\in I$. Hence there is some curve $\xi:I\to\gfrak_\alpha$ such that \cref{eq:recover-dynamics-symplectic:algebraic} holds. We must see that it is unique and smooth. First suppose there is some other curve $\eta:I\to\gfrak_\alpha$ fulfilling \cref{eq:recover-dynamics-symplectic:algebraic}. Then $\zeta(t):=\xi(t)-\eta(t)\in\ker\Phi^{\beta(t)}_*=\gfrak_{\beta(t)}$ for every $t\in I$. In the proof of \Cref{Group-Actions:injective-immersion} we saw that if $\zeta(t)_M(\beta(t))=0$, then $\exp(s\zeta(t))\in G_{\beta(t)}$ for all $s\in\R$. On the other hand, since $\zeta(t)_{\gfrak^*}(\alpha)=0$, because $\zeta(t)\in\gfrak_\alpha$, we also have that $\exp(s\zeta(t))\in G_\alpha$ for all $s\in\R$. Since the action of $G_\alpha$ on $\mu^{-1}(\alpha)$ is free, we obtain that $G_{\beta(t)}\cap G_\alpha=\{e\}$ for all $t\in I$, so that $\exp(s\zeta(t))=e$ for all $s\in\R$. This forces $\zeta(t)=0$ for all $t\in I$.

To see the smoothness, write $\xi(t)=\xi^i(t)\xi_i$, for $\{\xi_i\}$ a basis for $\gfrak_\alpha$ and some functions $\xi^i:I\to\R$. It is enough to see the smoothness of each $\xi^i$. We have that $\xi(t)_M(\beta(t))=\xi^i(t)(\xi_i)_M(\beta(t))$. Since for each $t$ the restriction of $\Phi^{\beta(t)}_*$ to $\gfrak_\alpha$ is injective, as has been shown in the previous paragraph, then the vectors $\{(\xi_i)_M(\beta(t))\}$ are linearly independent. In addition, $(\xi_i)_M(\beta(t))$ is the composition of the smooth maps $\beta:I\to M$ and $(\xi_i)_M:M\to TM$, so that $(\xi_i)_M(\beta(t))$ is smooth in $t$. Since $\xi^i(t)(\xi_i)_M(\beta(t))$ is also smooth in $t$, necessarily the functions $\xi^i$ are smooth, and this ends the proof.
\end{Proof}

The problem of finding the motion on the original system, therefore, reduces to solving the algebraic equation (\ref{eq:recover-dynamics-symplectic:algebraic}) and then the differential equation (\ref{eq:recover-dynamics-symplectic:differential}).

Recall that by \Cref{ex:symplectic-are-poisson} any symplectic manifold is a Poisson manifold with the Poisson bracket $\{f,g\}=-\omega(X_f,X_g)$. One could naturally ask how the Poisson brackets on $M$ and on $M\dbar_\alpha G$ are related. Let $f,h\in\Cinf(M\dbar_\alpha G)$ and let $F,H\in\Cinf(M)$ be $G$-invariant extensions of $\pi^*f$ and $\pi^*h$, respectively.

Using \Cref{thm:reduced-symplectic-dynamics}, we have that $X_F$ is $\pi$-related to $X_f$ and $X_H$ is $\pi$-related to $X_h$. Hence, it is immediate to see that if we write $\{\cdot,\cdot\}_\alpha$ for the Poisson bracket on $M\dbar_\alpha G$,
\[
\pi^*\{f,h\}_\alpha=i^*\{F,H\}.
\]


\section{Poisson Reduction}

\subsection{Marsden-Ratiu Reduction}

Poisson reduction as will be presented here was introduced by Marsden and Ratiu in 1986 in a short paper \cite{marsden-ratiu}. Of course, Poisson reduction should include the above Marsden-Weinstein reduction as a particular example. It gives, though, a more general construction, allowing for other interesting examples. 

Let $(M,\{\cdot,\cdot\})$ be a Poisson manifold with Poisson bivector field $\Pi$. Let $N\subseteq M$ be a submanifold with inclusion $i:N\to M$. We say that a subbundle $E\subseteq TM|_N$ is a \textbf{Poisson distribution} if 
\begin{enumerate}[label=(PD\arabic*),leftmargin=3\parindent]
	\item $E\cap TN$ is an integrable subbundle of $TN$, defining a foliation $\Fscr$ on $N$\footnote{See the Appendix for the definition of foliation.}, \label{poisson-distribution:1}
	\item the leaf space $N/\Fscr$ is a smooth manifold such that the projection $\pi:N\to N/\Fscr$ is a submersion, \label{poisson-distribution:2}
	\item for any $F,G\in\Cinf(M)$ such that $dF(E)=dG(E)=0$, we have that $d\bracket{F}{G}(E)=0$. \label{poisson-distribution:3}
\end{enumerate}

\begin{definition}
The triple $(M,N,E)$ is called \textbf{Poisson reducible} if there is a Poisson structure $\bracket{\cdot}{\cdot}_{N/\Fscr}$ on the leaf space $N/\Fscr$ such that for any $f,g\in\Cinf(N/\Fscr)$ and $F,G\in\Cinf(M)$ local extensions of $\pi^*f$ and $\pi^*g$, respectively, with $dF(E)=dG(E)=0$, we have that
\[
i^*\bracket{F}{G}=\pi^*\bracket{f}{g}_{N/\Fscr}. \qedhere
\]
\end{definition}

When we say that $F$ and $G$ are local extensions of $\pi^*f$ and $\pi^*g$ with $dF(E)=dG(E)=0$ we mean that this last condition must be satisfied locally but need not be satisfied throughout all of $N$. That is to say, when evaluating $\bracket{f}{g}_{N/\Fscr}$ at $\pi(p)\in N/\Fscr$, the extensions $F$ and $G$ must satisfy $dF(E)=dG(E)=0$ only in a neighborhood of $p$. Note that such local extensions always exist. Indeed, taking a slice chart for $N$ about $p\in N$ (that is, a chart $(U,\phi)$ such that $p\in U$ and $\phi(U\cap N)=\phi(U)\cap (\R^n\times\{0\})$, where $\dim N=n$), it suffices to see the existence of local extensions in the case where $M=\R^m$ for some $m$, $N=\R^n\times\{0\}$ and $E$ is a distribution of $k$-planes on $N$ ($k$ is the rank of $E$). If we denote by $\pi_1:\R^n\times\R^{m-n}\to\R^n$ the projection onto the first component and we are given a smooth function $f$ on $N$ constant on the leaves of the foliation induced by $E\cap TN$, then $F=\pi_1^*f$ is an extension of $f$ with $dF(E)=0$.

The reduction theorem states a necessary and sufficient condition for a triple $(M,N,E)$, where $E$ is a Poisson distribution over $N$, to be Poisson reducible. We denote by $E^\circ$ the annihilator of $E$ (defined pointwise as $(E^\circ)_p:=E_p^\circ\subseteq T_p^*M$).

\begin{theorem}[Marsden-Ratiu] \label{thm:marsden-ratiu}
The triple $(M,N,E)$ is Poisson reducible if and only if 
\begin{equation} \label{eq:poisson-reduction}
\Pi(E^\circ)\subseteq TN+E.
\end{equation}
\end{theorem}

\begin{Proof}
If $(M,N,E)$ is Poisson reducible, let $p\in N$, $\lambda\in E_p^\circ$ and $F\in\Cinf(M)$ such that $dF_p=\lambda$ and $dF(E)=0$. Let $\beta\in T_pN^\circ\cap E_p^\circ=(T_pN+E_p)^\circ$ and let $G\in\Cinf(M)$ be an extension of the zero function on $N$ such that $dG_p=\beta$ and $dG(E)=0$. Since $dF(E)=0$, we have that $F$ is constant on the leaves of $\Fscr$, and, hence, it descends to a smooth function $f\in\Cinf(N/\Fscr)$ such that $\pi^*f=F$ on a neighborhood of $p$. Then \cref{eq:poisson-reduction} follows from
\[
\Pi(\lambda)(\beta)=\Pi(dF,dG)(p)=i^*\bracket{F}{G}(p)=\pi^*\bracket{f}{0}_{N/\Fscr}(p)=0.
\]

Conversely, suppose that $\Pi(E^\circ)\subseteq TN+E$. Let $f,g\in\Cinf(N/\Fscr)$ and $F,G\in\Cinf(M)$ extensions of $\pi^*f$ and $\pi^*g$, respectively, such that $dF(E)=dG(E)=0$ about $p\in N$. Since $E$ is a Poisson distribution, by \ref{poisson-distribution:3} we know that $\bracket{F}{G}$ is constant along the leaves of $\Fscr$ and so gives rise to a well-defined function on $N/\Fscr$, $\bracket{f}{g}_{N/\Fscr}$, such that $\pi^*\bracket{f}{g}_{N/\Fscr}=i^*\bracket{F}{G}$. Define the bracket of $f$ and $g$ by $\bracket{f}{g}_{N/\Fscr}$.

To see that this definition does not depend on the choice of extensions, assume that $G'\in\Cinf(M)$ is another local extension of $\pi^*g$ with $dG'(E)=0$. Then $G-G'$ is constant on $N$ and therefore $d(G-G')|_N \in (TN+E)^\circ$. Thus, by \cref{eq:poisson-reduction},
\[
\Pi(dF|_N)(d(G-G')|_N)(p)=0,
\]
and $\bracket{\cdot}{\cdot}_{N/\Fscr}$ is well defined.

Bilinearity, skew-symmetry and the Leibniz rule for $\bracket{\cdot}{\cdot}_{N/\Fscr}$ are inherited directly from $\bracket{\cdot}{\cdot}$. For the Jacobi identity, we have already seen that $\bracket{F}{G}$ is a suitable local extension for $\pi^*\bracket{f}{g}_{N/\Fscr}$, so that for any $h\in\Cinf(N/\Fscr)$ and its extension $H\in\Cinf(M)$ with $dH(E)=0$,
\[
\pi^*\bracket{\bracket{f}{g}_{N/\Fscr}}{h}_{N/\Fscr}=i^*\bracket{\bracket{F}{G}}{H}.
\]
Hence, the Jacobi identity is also inherited from $\bracket{\cdot}{\cdot}$.
\end{Proof}

Also in this Poisson context, the book \cite{ortega-ratiu} explores what can be retained when the regularity conditions are dropped.

We present now some examples of Poisson reduction that are of interest.

\begin{example} \label{ex:quotient-by-canonical-action}
Let $(M,\bracket{\cdot}{\cdot})$ be a Poisson manifold and $G$ a Lie group acting canonically on $M$ (i.e., $\Phi_g$ preserves Poisson brackets for each $g\in G$). Let $N=M$ and $E_p:=T_p(G\cdot p)$ for each $p\in M$. The foliation $\Fscr$ is just the collection of $G$-orbits and $N/\Fscr=M/G$. Because the action is canonical, if $F,H\in\Cinf(M)$ are constant on the $G$-orbits then $\{F,H\}$ is so as well. Since $N=M$, obviously $\Pi(E^\circ)\subseteq TN+E$, so that $M/G$ is a Poisson manifold by \Cref{thm:marsden-ratiu}. For functions $f,g\in\Cinf(M/G)$, $\pi^*f$ and $\pi^*g$ are constant on the orbits, so that
\[
\pi^*\bracket{f}{g}_{M/G}=\bracket{\pi^*f}{\pi^*g}. \qedhere
\]
\end{example}

As with symplectic reduction, Poisson reduction also allows us to justify otherwise special Poisson structures on some objects. To see this, though, we first need the expression of the Poisson bracket on the cotangent bundle of a Lie group.

\begin{proposition} \label{prop:Poisson-bracket-cotangent-bundle-lie-group}
Let $G$ be a Lie group and $\gfrak$ its Lie algebra. Identify $T^*G$ with $G\times\gfrak^*$ via right translations. Let $(g,\alpha)\in G\times\gfrak^*$ and define $i_\alpha:G\to G\times\gfrak^*$ by $i_\alpha(h)=(h,\alpha)$ and $i_g:\gfrak^*\to G\times\gfrak^*$ by $i_g(\lambda)=(g,\lambda)$. Then if $F,H\in\Cinf(G\times\gfrak^*)$, we have that $dF_{(g,\alpha)}=d(i_\alpha^*F)_g+d(i_g^*F)_\alpha.$ Let $F_g:=d(i_\alpha^*F)_g$ and $F_\alpha:=d(i_g^*F)_\alpha$, and similarly for $H$. The canonical Poisson structure on $T^*G$ can be expressed as
\[
\bracket{F}{H}(g,\alpha)=H_g(R_{g*}F_\alpha)-F_g(R_{g*}H_\alpha)+\alpha([F_\alpha,H_\alpha]),
\]
where $F_\alpha$ and $H_\alpha$, which are linear functionals on $\gfrak^*$, are to be regarded as elements of $\gfrak$.
\end{proposition}

\begin{Proof}
Let $X_F(g,\alpha)=R_{g*}\xi+\beta$. Then, for any $R_{g*}\eta+\gamma\in T_gG\oplus\gfrak^*$,
\[
(i_{X_F}\omega)_{(g,\alpha)}(R_{g*}\eta+\gamma)=dF_{(g,\alpha)}(R_{g*}\eta+\gamma)=F_g(R_{g*}\eta)+F_\alpha(\gamma)=\omega(R_{g*}\xi+\beta,R_{g*}\eta+\gamma).
\]
Comparing with the expression of the symplectic form, \Cref{prop:symplectic-form-cotangent-bundle-lie-group}, we see that $\xi=F_\alpha$ and $\beta=-F_g\circ R_{g*}-\alpha\circ\ad(F_\alpha)$. Hence,
\[
X_F(g,\alpha)=R_{g*}F_\alpha-F_g\circ R_{g*}-\alpha\circ\ad(F_\alpha),
\]
so that
\[
\bracket{F}{H}(g,\alpha)=dH(X_F)(g,\alpha)=H_g(R_{g*}F_\alpha)-F_g(R_{g*}H_\alpha)-\alpha([F_\alpha,H_\alpha]),
\]
as wanted.
\end{Proof}

\begin{example}
Take \Cref{ex:quotient-by-canonical-action} with $M=T^*G\cong G\times\gfrak^*$ and right translation lifted to $T^*G$ as a $G$-action, which acts through
\[
(g,\alpha)\cdot h=(gh,\alpha).
\]
Since this action is symplectic, by \Cref{Moment-Maps:construction-corollary}, then, by \Cref{ex:symplectic-are-poisson}, the action is also canonical. Hence $T^*G/G\cong (G\times\gfrak^*)/G\cong \gfrak^*$ is a Poisson manifold. The projection $\pi:G\times\gfrak^*\to\gfrak^*$ is just projection onto the second factor. Then, if $F,H\in\Cinf(\gfrak^*)$, for any $\alpha\in\gfrak^*$ and $g\in G$ we have that
\[
\bracket{F}{H}_{\gfrak^*}(\alpha)=\bracket{\pi^*F}{\pi^*H}(g,\alpha).
\]
With the notation of \Cref{prop:Poisson-bracket-cotangent-bundle-lie-group} we have that $i_\alpha^*\pi^*F=(\pi\circ i_\alpha)^*F=F(\alpha)$ and $i_g^*\pi^*F=(\pi\circ i_g)^*F=F$, so that $(\pi^*F)_g=d(i_\alpha^*\pi^*F)_g=0$ and $(\pi^*F)_\alpha=d(i_g^*\pi^*F)_\alpha=dF_\alpha$, and similarly for $H$. Hence, we conclude that
\[
\bracket{F}{H}_{\gfrak^*}(\alpha)=-\alpha([dF_\alpha, dH_\alpha]).
\]
This Poisson structure on $\gfrak^*$ is known as the Lie-Poisson structure on $\gfrak^*$. It was first found by S. Lie \cite{lie-aleman} and then rediscovered by Berezin \cite{berezin}, both writing its expression in local coordinates. For a further account on the Lie-Poisson structure, see \cite{weinstein-poisson}.
\end{example}

\subsection{Reduced Poisson Dynamics}

As in the symplectic case, dynamics can also be reduced in the Poisson context. If $N\subseteq M$ is a submanifold, we say that a function $H\in\Cinf(M)$ preserves $N$ if $X_H(p)\in T_pN$ for all $p\in N$.

\begin{theorem}[Reduced Poisson dynamics] \label{thm:reduced-poisson-dynamics}
Let $(M,N,E)$ be a Poisson reducible triple and let $H\in\Cinf(M)$ be a function preserving $N$ and such that $dH(E)=0$. Then the flow of $X_H$ leaves $N$ invariant and induces a Hamiltonian flow on $N/\Fscr$ with Hamiltonian $h$ given by $i^*H=\pi^*h$, where $i:N\to M$ is the inclusion and $\pi:N\to N/\Fscr$ the projection. In particular, $X_H$ and $X_h$ are $\pi$-related. The function $h$ is called \textbf{reduced Hamiltonian} or energy.
\end{theorem}

\begin{Proof}
Since $N$ is preserved by $X_H$, clearly its flow, say $\Theta_t$, leaves $N$ invariant. Hence, it gives rise to a well-defined flow $\theta_t$ on $N/\Fscr$ such that $\pi\circ\Theta_t=\theta_t\circ\pi$. Let $X$ be the vector field associated to $\theta_t$. Then $X_H$ and $X$ are $\pi$-related: for any $p\in N$,
\[
X\circ\pi(p)=\derev{}{t}{0}\theta_t\circ\pi(p)=\derev{}{t}{0}\pi\circ\Theta_t(p)=\pi_*X_H(p).
\]
The fact that $H$ is constant on the leaves ensures the existence of a smooth function $h$ in $N/\Fscr$ uniquely determined by $i^*H=\pi^*h$. For any $f\in\Cinf(N/\Fscr)$ and $F\in\Cinf(M)$ its local extension with $dF(E)=0$, we have that
\begin{align*}
df(X)(\pi(p))&=X(\pi(p))f=(\pi_*X_H(p))f=df(\pi_*X_H(p))=d\pi^*f(X_H(p))\\
&=dF(X_H(p))=\bracket{H}{F}(p)=\bracket{h}{f}_{N/\Fscr}(\pi(p)).
\end{align*}
Thus, $X=X_h$, and we are done.
\end{Proof}

The reduced Poisson dynamics allows us to see that Marsden-Ratiu and Marsden-Weinstein reduction coincide in the symplectic case, in the following sense. Let $(M,\bracket{\cdot}{\cdot})$ be a Poisson manifold and $\mu:M\to\gfrak^*$ an equivariant moment map for some canonical action of a Lie group $G$ on $M$. Let $\alpha\in\gfrak^*$ be a clean value for $\mu$ and suppose that $G_\alpha$ acts freely and properly on $\mu^{-1}(\alpha)$. Then there is a unique Poisson structure $\bracket{\cdot}{\cdot}_\alpha$ on $M\dbar_\alpha G:=\mu^{-1}(\alpha)/G_\alpha$ such that if $f,h\in\Cinf(M\dbar_\alpha G)$ and $F,H\in\Cinf(M)$ are local extensions of $\pi^*f$ and $\pi^*h$, respectively, constant on the $G$-orbits on $\mu^{-1}(\alpha)$, meaning that $dF(T_p(G\cdot p))=dH(T_p(G\cdot p))=0$ for every $p\in\mu^{-1}(\alpha)$, then
\begin{equation} \label{eq:marsden-ratiu-relation-in-corollary}
\pi^*\bracket{f}{h}_\alpha=i^*\bracket{F}{H},
\end{equation}
where $\pi:\mu^{-1}(\alpha)\to M\dbar_\alpha G$ is the projection and $i:\mu^{-1}(\alpha)\to M$ the inclusion. If in addition $M$ is symplectic and $\bracket{\cdot}{\cdot}$ is the Poisson bracket associated to its symplectic structure, then $\bracket{\cdot}{\cdot}_\alpha$ is the Poisson bracket associated to the unique symplectic structure on $M\dbar_\alpha G$ given by the Marsden-Weinstein theorem (\Cref{Reduction:thm:Marsden-Weinstein}). For more on this, see {\cite[Ex. B]{marsden-ratiu}}.

\chapter{Implosion} \label{chap:implosion}

Symplectic implosion was introduced by Guillemin, Jeffrey and Sjamaar in 2002, \cite{guillemin-jeffrey-sjamaar}. Loosely speaking, it is a way of ``abelianizing'' the action of a compact Lie group on a symplectic manifold at the cost of introducing singularities in the manifold. That is, if $G$ is a compact Lie group acting in a Hamiltonian fashion on a symplectic manifold $(M,\omega)$, then the so-called imploded cross-section $\Mimpl$ inherits a Hamiltonian action of a maximal torus $T\subseteq G$ such that $M\dbar_\alpha G=\Mimpl\dbar_\alpha T$ for certain values of $\alpha$. The price to pay is that $\Mimpl$ will not be a symplectic manifold in general, but a stratified symplectic space.


\section{Symplectic Implosion}

We fix once and for all a compact connected Lie group $G$ and a maximal torus $T\subseteq G$. Let $\gfrak$ denote the Lie algebra of $G$ and $\tfrak\subseteq\gfrak$ the Lie algebra of $T$.

\subsection{Root Decomposition and Weyl Chambers}

We first review the root decomposition of $\gfrak$ and the concept of Weyl chambers. Here we mainly follow \cite{sepanski}. Let $\zfrak=\{\xi\in\gfrak:[\xi,\eta]=0,\text{ for all $\eta\in\gfrak$}\}$ be the center of $\gfrak$ and let $[\gfrak,\gfrak]$ be the ideal generated by $\{[\xi,\eta]:\xi,\eta\in\gfrak\}$. Then {\cite[Thm. 5.18]{sepanski}} tells us that $\gfrak=\zfrak\oplus[\gfrak,\gfrak]$ and that $[\gfrak,\gfrak]$ is semisimple. On the other hand, if we write $_\C$ for complexification, $\gfrak_\C$ always admits an $\Ad$-invariant (Hermitian) inner product (see {\cite[Lem. 5.6]{sepanski}}), so that the $\ad$ action is skew-Hermitian. Hence, $\ad\tfrak_\C$ is simultaneously diagonalizable and there is a finite set $\Delta\subseteq \tfrak_\C^*\smallsetminus\{0\}$ such that
\begin{equation}\label{eq:root-decomposition-compact}
\gfrak_\C=\tfrak_\C\oplus \bigoplus_{\alpha\in\Delta}\gfrak_\alpha,
\end{equation}
where $\gfrak_\alpha:=\{\xi\in\gfrak_\C:\ad\eta(\xi)=\alpha(\eta)\xi,\text{ for all $\eta\in\tfrak_\C$}\}\neq \emptyset$. \Cref{eq:root-decomposition-compact} is called the \textbf{root decomposition} of $\gfrak_\C$. The set $\Delta$ is called the set of \textbf{roots}, and, if we set $\tfrak':=\tfrak\cap[\gfrak,\gfrak]$, it spans $\tfrak'^*_\C$ ({\cite[Thm. 6.11]{sepanski}}).

Because of the skew-Hermiticity of $\ad\gfrak$, we have that $\alpha$ is imaginary-valued on $\tfrak$. Hence, by $\C$-linearity, any root $\alpha\in\Delta$ is completely determined by its restriction to $\tfrak$, and so we will interchangeably think of the roots as living in $\tfrak^*$ or $\tfrak_\C^*$. Using the root decomposition it is very easy to see that $\xi\in\tfrak$ lies in $\zfrak$ if and only if $\alpha(\xi)=0$ for all $\alpha\in\Delta$.

The \textbf{Killing form} is defined as the bilinear symmetric form on $\gfrak$ given by $B(\xi,\eta):=\tr(\ad\xi\circ\ad\eta)$. Restricted to $[\gfrak,\gfrak]$, it is negative definite (see {\cite[Thm. 6.16]{sepanski}}), and by $\C$-linearity it can be extended to a nondegenerate bilinear symmetric form on $[\gfrak_\C,\gfrak_\C]$. By taking any nondegenerate bilinear symmetric form on $\zfrak_\C$ (notice that any bilinear form on $\zfrak_\C$ is $\Ad$-invariant), we can extend $B$ further to a nondegenerate and $\Ad$-invariant bilinear symmetric form $\prodesc{\cdot}{\cdot}$ on all of $\gfrak_\C$ such that $\zfrak_\C$ and $[\gfrak_\C,\gfrak_\C]$ are orthogonal. For any $\lambda\in\gfrak^*_\C$, let $\zeta_\lambda\in\gfrak_\C$ be the unique element such that $\lambda(\xi)=\prodesc{\zeta_\lambda}{\xi}$ for all $\xi\in\gfrak_\C$; we can then define a nondegenerate and $\Ad^*$-invariant bilinear symmetric form on $\gfrak^*_\C$ as $\prodesc{\lambda}{\beta}:=\prodesc{\zeta_\lambda}{\zeta_\beta}$.

Observe that if $\lambda\in\tfrak_\C^*$, then $\zeta_\lambda\in\tfrak_\C$. Indeed, first of all note that in the root decomposition \eqref{eq:root-decomposition-compact}, the two terms in the direct sum are mutually orthogonal, since if $\xi\in\tfrak_\C$ and $\eta_\alpha\in\gfrak_\alpha$ for some $\alpha\in\Delta$, then for any $\zeta\in\tfrak$,
\[
0=\prodesc{\ad\zeta(\xi)}{\eta_\alpha}+\prodesc{\xi}{\ad\zeta(\eta_\alpha)}=\alpha(\zeta)\prodesc{\xi}{\eta_\alpha}.
\]
Hence, if we write $^\perp$ for the orthogonal space, $\lambda\in\tfrak_\C^*=\bigcap_{\alpha\in\Delta}\gfrak_\alpha^\circ$ if and only if $\zeta_\lambda\in\bigcap_{\alpha\in\Delta}\gfrak_\alpha^\perp=\tfrak_\C$.

We say that a subset $\Sigma\subseteq\Delta$ is a system of \textbf{simple roots} if it spans $\tfrak'^*$ (since $\alpha(\zfrak)=0$, it makes sense to think of the roots as elements of $\tfrak'^*$) and any root in $\Delta$ can be written as a linear combination of the elements of $\Sigma$ of either nonnegative or nonpositive coefficients. 
The set $\tfrak'^*\smallsetminus \bigcup_{\alpha\in\Delta}\alpha^\perp$ decomposes into some finite number of connected components, each of which is called an open \textbf{Weyl chamber}.

There is a bijection between systems of simple roots and Weyl chambers (see {\cite[Lem. 4.42]{sepanski}}). Given a system of simple roots $\Sigma$, we define the corresponding \textbf{fundamental closed Weyl chamber} as
\[
\tfrak_+^*:=\overline{\{\lambda\in\tfrak'^*: \prodesc{\lambda}{\alpha}>0,\text{ for all $\alpha\in\Sigma$.}\}}.
\]
We fix once and for all a system of simple roots $\Sigma$ and the corresponding fundamental closed Weyl chamber $\tfrak_+^*$.

\subsection{Symplectic Implosion}

Recall that for a group $H$, its commutator subgroup is defined as $[H,H]:=\langle ghg^{-1}h^{-1}: g,h\in H\rangle$. We give now the basic construction of symplectic implosion.

\begin{definition}
Let $(M,\omega)$ be a symplectic manifold with an equivariant moment map $\mu:M\to\gfrak^*$ arising from a Hamiltonian $G$-action $\Phi$. We define the following equivalence relation: $p,q\in \mu^{-1}(\tfrak^*_+)$ are equivalent if there is some $g\in [G_{\mu(p)},G_{\mu(p)}]$ such that $\Phi_g(p)=q$.

If we write $\sim$ for this equivalence relation, we define the \textbf{imploded cross-section} of $M$ as the quotient space $\Mimpl:=\mu^{-1}(\tfrak^*_+)/\sim$, equipped with the quotient topology.
\end{definition}

Indeed, $\sim$ is an equivalence relation because of the equivariance of $\mu$: if $p\sim q$ and $g\in[G_{\mu(p)},G_{\mu(p)}]$ is such that $\Phi_g(p)=q$, then $\mu(q)=\mu\circ\Phi_g(p)=\Ad^*(g)\mu(p)=\mu(p),$ so that the relation is symmetric.

If $\Sigma$ has $n$ roots, then a subset $\sigma\subseteq\tfrak_+^*$ is called a $(n-k)$-face if there is some subset $\Sigma_0\subseteq\Sigma$ with $k$ elements such that
\[
\sigma=\{\lambda\in\tfrak_+^*: \prodesc{\lambda}{\alpha}=0,\text{ for $\alpha\in\Sigma_0$}\}.
\]
We define the following partial order on the set of faces of $\tfrak_+^*$: if $\sigma$ and $\tau$ are faces of $\tfrak_+^*$ we say that $\sigma\leq\tau$ if $\sigma\subseteq \overline{\tau}$.

To give a first description of $\Mimpl$, we need the fact that the isotropy group with respect to the coadjoint action is constant on the faces of the fundamental Weyl chamber.

\begin{proposition} \label{prop:isotropy-groups-faces}
The isotropy group by the coadjoint action is the same for all the elements of a given face in $\tfrak_+^*$. Moreover, $T$ is contained in any of these isotropy groups and it is the isotropy group of any element in the interior face of $\tfrak_+^*$.
\end{proposition}

\begin{Proof}
Since $G$ is connected, then the isotropy groups for the coadjoint action are also connected (see the first remark after {\cite[Lem. 2.3.2]{guillemin-lerman-sternberg}}). Hence, the equality of isotropy groups is equivalent to the equality of isotropy algebras. 

Fix a $(n-k)$-face $\sigma\subseteq\tfrak_+^*$ defined by a subset $\Sigma_0\subseteq\Sigma$. Let $\lambda\in\sigma$. By \Cref{Moment-Maps:interesting-properties} and \Cref{Group-Actions:isotropy-algebra:p-tangent-orbit}, the complexification of its isotropy algebra is
\[
(T_eG_\lambda)_\C=\{\xi\in\gfrak_\C: \lambda\circ\ad\xi=0\},
\]
where here we are thinking of $\lambda$ as an element in the annihilator of $\zfrak_\C\oplus\bigoplus_{\alpha\in\Delta} \gfrak_\alpha$. If we define $\Delta_0:=\{\alpha\in\Delta: B(\lambda,\alpha)=0\}$, we will show that
\[
(T_eG_\lambda)_\C= \tfrak_\C\oplus\bigoplus_{\alpha\in\Delta_0}\gfrak_\alpha,
\]
and this will give the results, since $\Delta_0$ is totally determined by $\Sigma_0$, which is independent of $\lambda$.

If $\xi\in\gfrak_\C$, then $\lambda\circ\ad\xi=0$ if and only if for all $\eta\in\gfrak_\C$ we have that
\[
\lambda([\xi,\eta])=\prodesc{\zeta_\lambda}{[\xi,\eta]}=-\prodesc{[\xi,\zeta_\lambda]}{\eta}=0,
\]
i.e., if and only if $\ad\zeta_\lambda(\xi)=0$. Hence, in fact $(T_eG_\lambda)_\C=\ker\ad\zeta_\lambda$. 

For any $\xi=\xi_0+\sum_{\alpha\in\Delta}\xi_\alpha$, where $\xi_0\in\tfrak_\C$ and $\xi_\alpha\in\gfrak_\alpha$, we have that $[\zeta_\lambda,\xi]=\sum_{\alpha\in\Delta}\alpha(\zeta_\lambda)\xi_\alpha$, because $\zeta_\lambda\in\tfrak_\C$. Therefore, $\xi\in\ker\ad\zeta_\lambda$ if and only if $\xi=\xi_0+\sum_{\alpha\in\Delta_0}\xi_\alpha$, because $\alpha(\zeta_\lambda)=B(\alpha,\lambda)$, and this ends the proof.
\end{Proof}

Hence, for a face $\sigma$, we may write without confusion $G_\sigma$ to refer to the isotropy group of any of its elements. If we call $F$ the set of faces of $\tfrak_+^*$, then $\tfrak_+^*=\bigsqcup_{\sigma\in F}\sigma$. Thus, we may rewrite the imploded cross-section (set-theoretically) as 
\begin{equation}\label{eq:implosion-set-theoretic}
\Mimpl=\bigsqcup_{\sigma\in F} \mu^{-1}(\sigma)/[G_\sigma,G_\sigma].
\end{equation}


\section{Cross-section Theorems and Symplectic Decomposition}\label{sect:cross-section-theorems}

We will now see that each of the disjoint sets of \cref{eq:implosion-set-theoretic} is actually a (possibly singular) symplectic quotient. The basic tool to see this is Guillemin and Sternberg's cross-section theorem, {\cite[Thm. 26.7]{guillemin-sternberg}}, or rather an improved version due to Lerman, Meinrenken, Tolman and Woodward, {\cite[Thm. 3.8]{lerman-meinrenken-tolman-woodward}}. It is a way of obtaining symplectic submanifolds through the moment map.

Let $(M,\omega)$ be a symplectic manifold with a Hamiltonian $G$-action, where $G$ is Lie group with Lie algebra $\gfrak$. Let $\mu:M\to \gfrak^*$ be the associated equivariant moment map. For $\lambda\in\gfrak^*$ we write $G\cdot\lambda$ for the coadjoint orbit.

\begin{theorem}[Guillemin-Sternberg] \label{thm:guillemin-sternberg-cross-section}
Let $\lambda\in\gfrak^*$ and $Z\subseteq \gfrak^*$ a submanifold perfectly transverse to $G\cdot \lambda$ at $\lambda$, that is, such that $T_\lambda Z \oplus T_\lambda (G\cdot \lambda)=\gfrak^*$. Then, if $p\in M$ with $\mu(p)=\lambda$, there is a neighborhood $U$ of $p$ such that $\mu^{-1}(Z)\cap U$ is a symplectic submanifold of $M$.
\end{theorem}

\begin{Proof}
We follow the proof in \cite{guillemin-sternberg}, correcting a mistake at the end of the proof. We repeatedly use \Cref{Moment-Maps:interesting-properties}, \Cref{Group-Actions:isotropy-algebra:p-tangent-orbit} and \Cref{lemma:kermu:symplectic}.

First of all, since $T_\lambda(G\cdot \lambda)=\mu_*(T_p(G\cdot p))$, $\mu$ is transverse to $Z$. Hence in a neighborhood of $p$, we have that $\mu^{-1}(Z)$ is a submanifold and $T_p\mu^{-1}(Z)=\mu_*^{-1}(T_\lambda Z)$. To see that in fact it is a symplectic manifold it suffices to see that $W_p:=T_p\mu^{-1}(Z)$ is a symplectic subspace, i.e., $W_p\cap W_p^\perp=0$.

Since $\ker \mu_*(p)=\gfrak_M(p)^\perp$, then $\gfrak_M(p)^\perp\subseteq \mu_*^{-1}(T_\lambda Z)=W_p$. Since $T_\lambda Z \cap T_\lambda(G\cdot \lambda)=0$, we also have that
\begin{align*}
W_p\cap \gfrak_M(p)&=\{\xi_M(p):\xi\in\gfrak, \mu_*\xi_M(p)\in T_\lambda Z\} =\{\xi_M(p):\xi\in\gfrak, \mu_*\xi_M(p)=0\}\\
& = \gfrak_M(p)\cap \ker\mu_*(p)=\gfrak_M(p)\cap \gfrak_M(p)^\perp.
\end{align*}
Therefore,
\[
W_p\cap W_p^\perp=W_p\cap (W_p^\perp\cap \gfrak_M(p))=\gfrak_M(p)\cap W_p^\perp \cap \gfrak_M(p)^\perp. 
\]

If $\xi\in\gfrak$ is such that $\xi_M(p)\in W_p^\perp$, then for all $v\in T_\lambda Z\cap \im\mu_*(p)$, that is, such that $v=\mu_*w$ for $w\in W_p$, we have that
\[
\xi(v)=\mu_*w(\xi)=d\muhat(\xi)(w)=\omega(\xi_M(p),w)=0,
\]
and hence $\xi\in(T_\lambda Z\cap \im\mu_*(p))^\circ$. On the other hand, if $\xi_M(p)\in \gfrak_M(p)^\perp=\ker \mu_*(p)$, then $\xi_{\gfrak^*}(\lambda)=0$, so that $\xi\in\gfrak_\lambda=(T_\lambda(G\cdot \lambda))^\circ$. 

Since 
\[
\im\mu_*(p)=T_\lambda Z\cap \im\mu_*(p)+T_\lambda (G\cdot \lambda),
\] 
we conclude that $\xi\in(\im\mu_*(p))^\circ=\gfrak_p$, i.e., $\xi_M(p)=0$, and this completes the proof.
\end{Proof}

If $G$ is compact and connected, which from now on we suppose, for any $\lambda\in\gfrak^*$ there is always a submanifold $S\subseteq\gfrak^*$ perfectly transverse to $\lambda$, called the natural slice. What the refinement of the cross-section theorem provided by Lerman, Meinrenken, Tolman and Woodward affirms is that $\mu^{-1}(S)$ is a symplectic submanifold not only locally, but also globally.

Before giving the definition of a slice to an action, recall that if $\pi:P\to M$ is a principal $G$-bundle and $N$ is a manifold with a left $G$-action, then the associated bundle with fiber $N$ is $P\times_G N:=(P\times N)/G$, where the $G$-action on $P\times N$ is $(p,n)\cdot g:=(p\cdot g,g^{-1}\cdot n)$ and the projection is given by $\pi'([p,n]):=\pi(p)$.\footnote{See the Appendix for the definition of principal $G$-bundle and, for instance, {\cite[Sect. I.5]{kobayashi-nomizu}} for the proof that the associated bundle is a fiber bundle.}

\begin{definition}
Let $M$ be a manifold with a $G$-action. If $p\in M$, we say that a submanifold $S\subseteq M$ is a \textbf{slice} for the $G$-action at $p$ if 
\begin{enumerate}[label=(S\arabic*),leftmargin=3\parindent]
	\item $S$ is $G_p$-invariant, \label{Implosion:def:slice:1}
	\item $G\times_{G_p}S\to M$ given by $[g,s]\mapsto \Phi_g(s)$ is a diffeomorphism onto its image, and \label{Implosion:def:slice:2}
	\item $G\cdot S$ is an open neighborhood for $G\cdot p$. \label{Implosion:def:slice:3} \qedhere
\end{enumerate}
\end{definition}

Notice that \ref{Implosion:def:slice:2} makes sense because the action of $G_p$ on $G$ by right translations defines a principal bundle $G\to G/G_p$ and \ref{Implosion:def:slice:1} ensures the existence of a $G_p$-action on $S$. Note also that by \ref{Implosion:def:slice:2}, for any $s\in S$, if $\Phi_g(s)\in S$ for some $g\in G$, then $g\in G_p$. Therefore $(G\cdot s)\cap S=G_p\cdot s$ and $G_s\subseteq G_p$.

The way one must think of a slice at $p$ is that it is somehow transverse to the orbits near $p$, and so it can be used to ``parameterize'' them in a neighborhood of $G\cdot p$ via the diffeomorphism $G\times_{G_p}S\to G\cdot S$, being this parameterization ``degenerate'' by $G_p$, in the sense that
two points in $S$ represent the same orbit if and only if there is some element of $G_p$ taking one to the other.

The slice is perfectly transverse to $G\cdot p$ at $p$, as the following concatenation of identifications shows
\begin{align*}
T_pM&=T_p(G\cdot S)=T_{[e,p]}(G\times_{G_p}S)\\
&=(T_eG\oplus T_pS)/(T_{(e,p)}(G_p\cdot (e,p))) \\
&=(T_eG\oplus T_pS)/(T_eG_p\oplus \{*\}) \\
&=T_eG/T_eG_p\oplus T_pS=T_p(G\cdot p)\oplus T_pS.
\end{align*}

It is the case that if $S$ is a $G_p$-invariant submanifold perfectly transverse to $G\cdot p$ at $p$ fulfilling that if $s\in S$ and $\Phi_g(s)\in S$ for some $g\in G$, then $g\in G_p$, this is sufficient for $S$ to be a slice at $p$ (see {\cite[Sect. 2.3.2]{guillemin-lerman-sternberg}}).
This allows us to construct a slice for the coadjoint action at every point $\lambda\in\tfrak_+^*$, called the \textbf{natural slice} at $\lambda$.

\begin{proposition}[\protect{{\cite[Lem. 2.3.2]{guillemin-lerman-sternberg}}}]
Let $\lambda$ be an element in a face $\sigma$ of $\tfrak_+^*$. The set 
\[
S_\sigma:=G_\sigma\cdot\{\beta\in\tfrak_+^*: G_\beta\subseteq G_\sigma\}
\]
is a slice for the coadjoint action at $\lambda$. 
Moreover,
\[
S_\sigma=G_\sigma\cdot \bigcup_{\sigma\leq \tau}\tau.
\]
\end{proposition}

We are now ready to state a second version of the cross-section theorem. 

\begin{theorem}[\protect{{\cite[Thm. 3.8]{lerman-meinrenken-tolman-woodward}}}] \label{thm:lmtw-cross-section}
Let $\sigma$ be a face of $\tfrak_+^*$, and let $S_\sigma$ be the natural slice for the coadjoint action at any point in $\sigma$. Then the \textbf{symplectic cross-section} $M_\sigma:=\mu^{-1}(S_\sigma)$ is a $G_\sigma$-invariant symplectic submanifold of $M$ and the restriction of $\mu$ to $M_\sigma$ is a moment map for the $G_\sigma$-action on it.
\end{theorem}

Let $\sigma$ be a face in $\tfrak_+^*$. Since $G_\sigma$ is compact, then $\gfrak_\sigma=\zfrak_\sigma\oplus[\gfrak_\sigma,\gfrak_\sigma]$, where $\zfrak_\sigma$ is the center of $\gfrak_\sigma$. Hence $\gfrak_\sigma^*=\zfrak_\sigma^*\oplus[\gfrak_\sigma,\gfrak_\sigma]^*$. If we let $\mu_\sigma:=\mu|_{M_\sigma}$, then the composition of $\mu_\sigma$ with the projection $\gfrak_\sigma^*\to[\gfrak_\sigma,\gfrak_\sigma]^*$ is a moment map for the action of $[G_\sigma,G_\sigma]$ on $M_\sigma$ whose zero locus is
\[
\mu_\sigma^{-1}(\zfrak_\sigma^*)=\mu^{-1}(\zfrak_\sigma^*)\cap M_\sigma=\mu^{-1}(\zfrak_\sigma^*\cap S_\sigma).
\]
Following {\cite[Rmk. 3.9c]{lerman-meinrenken-tolman-woodward}}, it can be seen that $\zfrak_\sigma^*$ is in fact the linear span of $\sigma$, so that $\zfrak_\sigma^*\cap S_\sigma=\sigma$, and, hence, \cref{eq:implosion-set-theoretic} can be rewritten as
\begin{equation}\label{eq:implosion-symp-quotients}
\Mimpl=\bigsqcup_{\sigma\in F}\mu^{-1}_\sigma(\zfrak_\sigma^*)/[G_\sigma,G_\sigma]=\bigsqcup_{\sigma\in F} M_\sigma\dbar_0[G_\sigma,G_\sigma].
\end{equation}

Not all the quotients in \cref{eq:implosion-symp-quotients} will be symplectic manifolds, some of them will be singular. In order to obtain a decomposition into symplectic manifolds one should further decompose each cross-section $M_\sigma$ into submanifolds by orbit type and then do the reduction. For this, one uses the following theorem, which is generalization of the Marsden-Weinstein reduction theorem.

\begin{theorem}[\protect{{\cite[Thm. 2.1]{sjamaar-lerman}}}] \label{thm:sjamaar-lerman}
Let $(M,\omega)$ be a symplectic manifold with an equivariant moment map $\mu:M\to\gfrak^*$ arising from a Hamiltonian $G$-action, where $G$ is compact. For any closed subgroup $H\leq G$ define the \textbf{orbit type $\boldsymbol{(H)}$ stratum} 
\[
M_{(H)}:=\{p\in M: G_p\text{ is conjugate to $H$}\}.
\]
Then $M_{(H)}\cap\mu^{-1}(0)$ is a submanifold such that 
\[
(M_{(H)}\cap\mu^{-1}(0))/G
\]
is a symplectic manifold with symplectic form $\omega_{(H)}$ fulfilling $\pi^*\omega_{(H)}=i^*\omega$, where $\pi:M_{(H)}\cap\mu^{-1}(0)\to (M_{(H)}\cap\mu^{-1}(0))/G$ is the projection and $i:M_{(H)}\cap\mu^{-1}(0)\to M$ the inclusion. Hence, the reduced space $M\dbar_0G$ decomposes into a disjoint union of symplectic manifolds
\[
M\dbar_0G=\bigsqcup_{H\leq G} (M_{(H)}\cap\mu^{-1}(0))/G.
\]
\end{theorem}

Therefore, for any closed subgroup $H$ of $[G_\sigma,G_\sigma]$ we have that $M_{\sigma,(H)}\cap\mu^{-1}(\sigma)$ is a submanifold such that $\Mimpl^{\sigma,(H)}:=(M_{\sigma,(H)}\cap\mu^{-1}(\sigma))/[G_\sigma,G_\sigma]$ is a symplectic manifold. Then we can partition the decomposition \eqref{eq:implosion-symp-quotients} further into symplectic manifolds
\begin{equation}\label{eq:implosion-symp-manifolds}
\Mimpl=\bigsqcup_{\substack{\sigma\in F\\ H\leq [G_\sigma,G_\sigma]}}\Mimpl^{\sigma,(H)}.
\end{equation}
In fact, this decomposition turns out to be a stratification for $\Mimpl$, see {\cite[Sect. 5]{guillemin-jeffrey-sjamaar}}.


\section{Abelianization}

We said that symplectic implosion is a way of abelianizing the action at the cost of introducing singularities in the manifold. We have already given a precise meaning to the ``introducing singularities'' part, and we would like now to turn to the ``abelianizing'' part. Here we essentially follow {\cite[Sect. 3]{guillemin-jeffrey-sjamaar}}. 

Let $\sigma$ be a face in $\tfrak_+^*$ and $\lambda\in\sigma$. To see that $M\dbar_\lambda G= \Mimpl\dbar_\lambda T$, let us define the Hamiltonian $T$-action on $\Mimpl$ in the following manner: first of all, if $p,q\in\mu^{-1}(\tfrak_+^*)$ are such that $p\sim q$, then $\Phi_g(p)\sim \Phi_g(q)$ for any $g\in T$, since $T$ normalizes $[G_{\mu(p)},G_{\mu(p)}]$, because $T\subseteq G_{\mu(p)}$, by \Cref{prop:isotropy-groups-faces}. Hence the $T$-action on $\mu^{-1}(\tfrak_+^*)$ descends to a $T$-action on $\Mimpl$ that is smooth on each of the symplectic manifolds into which $\Mimpl$ decomposes. 

If we call $\pi:\mu^{-1}(\tfrak_+^*)\to \Mimpl$ the projection, then $\mu$ induces a smooth map $\muimpl:\Mimpl\to\tfrak_+^*$ such that $\muimpl\circ\pi=\mu$. This $\muimpl$ is a moment map for the $T$-action on $\Mimpl$. To see this, consider the action of $G_\sigma$ on $\mu^{-1}(\sigma)$; it descends to an action of $G_\sigma/[G_\sigma,G_\sigma]$ on $M_\sigma\dbar_0[G_\sigma,G_\sigma]$, and since there is a natural projection $T\to G_\sigma/[G_\sigma,G_\sigma]$, because $T\subseteq G_\sigma$, this induces an action of $T$ on $M_\sigma\dbar_0[G_\sigma,G_\sigma]$. Once again because $T$ normalizes $[G_\sigma, G_\sigma]$, it is easy to see that $T$ preserves $M_{\sigma,(H)}$, so that $T$ acts in fact on each $\Mimpl^{\sigma,(H)}\subseteq M_\sigma\dbar_0[G_\sigma,G_\sigma]$. Since $\muimpl\circ\pi=\mu$, it is immediate that $\muhatimpl(\xi)\circ\pi=\muhat(\xi)$ for each $\xi\in\tfrak'$, and hence $\muimpl$ is a moment map for the $T$ action on $\Mimpl^{\sigma,(H)}$. Symplectic reduction of $\Mimpl$ by the $T$-action along $\lambda$ is defined by doing reduction on each $\Mimpl^{\sigma,(H)}$.

For simplicity, we will assume that $\lambda$ is such that all the points in $\mu^{-1}(\lambda)$ are of the same orbit type with respect to the $G_\sigma$-action, so that by \Cref{thm:sjamaar-lerman}, $M\dbar_\lambda G$ is a smooth manifold. 

\begin{theorem}[\protect{{\cite[Thm. 3.4]{guillemin-jeffrey-sjamaar}}}]
Let $\sigma$ be a face in $\tfrak_+^*$ and $\lambda\in\sigma$ such that all the points in $\mu^{-1}(\lambda)$ are of the same orbit type with respect to the $G_\sigma$-action. Then $M\dbar_\lambda G$ and $\Mimpl\dbar_\lambda T$ are symplectomorphic.
\end{theorem}

\begin{SProof}
By reduction by stages (see \cite{reduction-by-stages} for a painstaking analysis of symplectic reduction by stages), $M\dbar_\lambda G=\mu^{-1}(\lambda)/G_\sigma$ is symplectomorphic to
\begin{equation}\label{eq:reduction-by-stages}
(\mu^{-1}(\lambda)/[G_\sigma,G_\sigma])/T_\lambda,
\end{equation}
where $T_\lambda$ is the isotropy group of $\lambda$ with respect to the $\Ad^*(T)$-action on $\tfrak_+^*$.
Since all the elements in $\mu^{-1}(\lambda)$ are of the same orbit type, by \Cref{thm:sjamaar-lerman}, we have that $\muimpl^{-1}(\lambda)=\mu^{-1}(\lambda)/[G_\sigma,G_\sigma]$ is a symplectic manifold, whose symplectic reduction by $T$ along $\lambda$ is precisely \cref{eq:reduction-by-stages}. Since $\Mimpl\dbar_\lambda T=\muimpl^{-1}(\lambda)/T_\lambda$, we have the result.
\end{SProof}

In the general case where the orbit type of the points in $\mu^{-1}(\lambda)$ is not constant, one does the same reasoning on each of the disjoint components $\mu^{-1}(\lambda)\cap M_{\sigma,(H)}$. 


\section{Poisson Transversals} \label{sect:poisson-transversals}

We would now like to see what of the implosion construction survives in the Poisson setting. Since the Guillemin-Sternberg cross-section theorem was somehow the origin of it all, in the sense that it allowed for the decomposition of $\Mimpl$ into symplectic manifolds, we focus our attention mainly on generalizing this theorem. Since the Guillemin-Sternberg theorem concerned symplectic submanifolds, we first need to establish which kind of submanifolds, of the many different types that can be defined on Poisson manifolds, will be the correct one. 

Let $(M,\Pi)$ be a Poisson manifold with a Hamiltonian $G$-action, where $G$ is a Lie group with Lie algebra $\gfrak$. Let $\mu:M\to \gfrak^*$ be the associated equivariant moment map. 

\begin{definition}
A submanifold $N\subseteq M$ is called a \textbf{Poisson submanifold} if $\Pi$ is tangent to $N$, in the sense that $\Pi_p\in\Lambda^2T_pN$ for each $p\in N$. Then $\Pi|_N$ gives rise to a Poisson structure $\Pi_N$ on $N$ such that $\Pi_N$ and $\Pi$ are $i$-related, where $i:N\to M$ is the inclusion, meaning that $i_*(\Pi_N(p))=\Pi(p)$ for all $p\in N$. 
\end{definition}

Poisson submanifolds are somehow the most natural submanifolds to consider in a Poisson manifold. More natural in the sense that it obviously inherits a Poisson structure from the ambient space. Observe that a submanifold $N$ is Poisson if and only if $X_f$ is tangent to $N$ for every $f\in\Cinf(M)$, because $\Pi(T^*_pM|_N)=\operatorname{span}\{X_f(p):f\in\Cinf(M)\}$ for any $p\in N$.

The most naive approach to a generalization of the Guillemin-Sternberg cross-section theorem to the Poisson case would say that if we have $\lambda\in\gfrak^*$ and $Z\subseteq\gfrak^*$ a submanifold perfectly transverse to $G\cdot\lambda$ at $\lambda$, then $\mu^{-1}(Z)$ is a Poisson submanifold locally around some $p$ with $\mu(p)=\lambda$. But this is certainly not true, because if $\xi\in\gfrak$, then $X_{\muhat(\xi)}$ cannot be tangent to $\mu^{-1}(Z)$ at $p$, because $T_\lambda\mu^{-1}(Z)=\mu_*^{-1}(T_\lambda Z)$ and 
\[
\mu_*X_{\muhat(\xi)}(p)=\mu_*\xi_M(p)=\xi_{\gfrak^*}(\lambda)\in T_\lambda(G\cdot\lambda),
\]
and $T_\lambda(G\cdot\lambda)$ is precisely the direct sum complement of $T_\lambda Z$! 

Thus, the generalization must be done with greater care. As we will see, the correct notion is that of a Poisson tranversal (sometimes called cosymplectic submanifolds).

\begin{definition}
A submanifold $N\subseteq M$ is called a \textbf{Poisson transversal} if 
\[
TM|_N=TN+\Pi(TN^\circ). \qedhere
\]
\end{definition}

\begin{example}
Consider the Poisson manifold $(\R^5, \partial_1\wedge\partial_2+\partial_3\wedge\partial_4)$. The Poisson submanifolds are just the hyperplanes $\{x^5=k\}$, for fixed $k\in\R$. There are several Poisson transversals: the lines $\{(x^1,x^2,x^3,x^4)=(k^1,k^2,k^3,k^4)\}$ for fixed $k^1,k^2,k^3,k^4\in\R$, and the 3-planes $\{(x^1,x^2)=(k^1,k^2)\}$, for fixed $k^1,k^2\in\R$, and $\{(x^3,x^4)=(k^3,k^4)\}$, for fixed $k^3,k^4\in\R$, to name a few. 
\end{example}

There are some equivalent ways to characterize this type of submanifolds.

\begin{proposition}\label{PoissonTransversals:charac}
Let $N\subseteq M$ be a submanifold. Then the following are equivalent:
\begin{enumerate}
	\item $N$ is a Poisson transversal,
	\item $TM|_N=TN\oplus\Pi(TN^\circ)$, \label{PoissonTransversals:charac:2}
	\item $TN^\circ\cap \Pi^{-1}(TN)=0$, \label{PoissonTransversals:charac:3}
	\item $T^*M|_N=TN^\circ \oplus \Pi^{-1}(TN)$. \label{PoissonTransversals:charac:4}
\end{enumerate}
\end{proposition}

\begin{Proof}
If $N$ is a Poisson transversal and $\Pi(\lambda)\in T_pN\cap \Pi(T_pN^\circ)$ for some $p\in N$ and $\lambda\in T_pN^\circ$, then for every $w\in T_pM$ there are $v\in T_pN$ and $\alpha\in T_pN^\circ$ such that $w=v+\Pi(\alpha)$, and, hence,
\[
\lambda(w)=\lambda(\Pi(\alpha))=-\alpha(\Pi(\lambda))=0,
\]
so that $\lambda=0$. This shows the equivalence with \cref{PoissonTransversals:charac:2}. 

\Cref{PoissonTransversals:charac:3,PoissonTransversals:charac:4} are obtained taking annihilators and dualizing \ref{PoissonTransversals:charac:2}, respectively, since $\Pi(TN^\circ)^\circ=\Pi^{-1}(TN)$.
\end{Proof}

Though less obviously than in the case of a Poisson submanifold, Poisson transversals also carry a natural Poisson structure inherited from the ambient manifold. Because of the fact that the sum of $TN$ and $\Pi(TN^\circ)$ in a Poisson transversal is actually direct,
we can identify $T^*N$ with $\Pi^{-1}(TN)$ and $\Pi(TN^\circ)^*$ with $TN^\circ$, and therefore split $\Pi|_N$ in a suitable sense.

\begin{lemma} \label{PoissonTransversals:splitting-of-pi}
Let $N\subseteq M$ be a Poisson transversal. There are $\Pi_N\in\Xfrak^2(N)$ and $V\in\Gamma(\Lambda^2\Pi(TN^\circ))$ such that 
\[
\Pi|_N=\Pi_N+V,
\]
in the sense that for every $p\in N$ and every $\lambda=\lambda_1+\lambda_2$ and $\alpha=\alpha_1+\alpha_2$ in $T_p^*M$ with $\lambda_1,\alpha_1\in \Pi^{-1}(T_pN)$ and $\lambda_2,\alpha_2\in T_pN^\circ$,
\[
\Pi(\lambda,\alpha)=\Pi_N(\lambda_1,\alpha_1)+V(\lambda_2,\alpha_2).
\]
\end{lemma}

\begin{Proof}
With the identification $T^*N\cong \Pi^{-1}(TN)$, we define $\Pi_N(\lambda_1,\alpha_1):=\Pi(\lambda_1,\alpha_1)$, where on the left-hand side $\lambda_1$ and $\alpha_1$ are viewed as elements in $T^*N$ and on the right-hand side as elements in $\Pi^{-1}(TN)$. Similarly, we define $V(\lambda_2,\alpha_2):=\Pi(\lambda_2,\alpha_2)$. Then, since
\[
\Pi(\lambda_1,\alpha_2)=\alpha_2(\Pi(\lambda_1))=0
\]
because $\alpha_2\in T_pN^\circ$ and $\lambda_1\in\Pi^{-1}(T_pN)$, and $\Pi(\lambda_2,\alpha_1)=0$ as well, we obtain
\[
\Pi(\lambda,\alpha)=\Pi(\lambda_1,\alpha_1)+\Pi(\lambda_2,\alpha_2)=\Pi_N(\lambda_1,\alpha_1)+V(\lambda_2,\alpha_2),
\]
as wanted. $\Pi_N$ and $V$ are obviously smooth, since they are the composition of a bundle isomorphism with the bivector field $\Pi$.
\end{Proof}

We will now show that in fact $\Pi_N$ is a Poisson bivector field, defining the sought Poisson structure on $N$. Before proceeding to the proof, we need a technical lemma, whose proof is immediate because the Schouten-Nijenhuis bracket is defined pointwise.

\begin{lemma} \label{PoissonTransversals:tech-lemma}
Let $N\subseteq M$ be a submanifold and let $X\in\Xfrak^n(M)$ and $Y\in\Xfrak^m(M)$ be multivector fields tangent to $N$, in the sense that $X|_N\in \Xfrak^n(N)$ and $Y|_N\in\Xfrak^m(N)$. Then $[X,Y]|_N=[X|_N,Y|_N]$.
\end{lemma}

\begin{theorem}
Let $N\subseteq M$ be a Poisson transversal. The bivector field $\Pi_N$ defines a Poisson structure on $N$.
\end{theorem}

\begin{Proof}
It only remains to show that $[\Pi_N,\Pi_N]=0$. We follow the proof in \cite{xu}. Let $V$ be as in \Cref{PoissonTransversals:splitting-of-pi}, and let $V=\sum_i X_i\wedge Y_i$ for some $X_i,Y_i\in\Gamma(\Pi(TN^\circ))$ (we can do this at least locally). Let $\Xline_i,\Yline_i\in\Xfrak(M)$ be local extensions of $X_i, Y_i$, and let $\Vline=\sum_i\Xline_i\wedge \Yline_i$. If $\Piline_N=\Pi-\Vline$, then $\Vline|_N=V$ and $\Piline_N|_N=\Pi_N$. The fact that $[\Pi,\Pi]=0$ implies that
\[
[\Piline_N,\Piline_N]=-2[\Pi,\Vline]+[\Vline,\Vline]=-2[\Piline_N,\Vline]-[\Vline,\Vline].
\]
If $p:TM|_N\to TN$ denotes the projection given by the decomposition $TM|_N=TN\oplus \Pi(TN^\circ)$, as well as its extension to multivector fields, then it is clear that, by \Cref{PoissonTransversals:tech-lemma},
\[
p([\Piline_N,\Vline]|_N)=p([\Pi_N,V])=0 \quad\text{and}\quad p([\Vline,\Vline]|_N)=p([V,V])=0.
\]
Hence we have that $p([\Piline_N,\Piline_N]|_N)=0$. We conclude that
\[
[\Pi_N,\Pi_N]=p([\Pi_N,\Pi_N])=p([\Piline_N,\Piline_N]|_N)=0. \qedhere
\]
\end{Proof}

In terms of smooth functions on $N$, the Poisson structure can be described in the following way: let $f,g\in\Cinf(N)$ and let $F,G\in\Cinf(M)$ be extensions of $f$ and $g$, respectively, such that $dF$ or $dG$ annihilates $\Pi(TN^\circ)$ locally, then
\[
\{f,g\}_N=\{F,G\}.
\]
In the case $(M,\Pi)$ were symplectic, with symplectic form $\omega$, then, since $TM|_N=TN\oplus\Pi(TN^\circ)=TN\oplus\omega^{-1}(TN^\circ)=TN\oplus TN^\perp$, we have that $N$ is a symplectic submanifold. We now see that Poisson transversals are the correct candidate for a Poisson cross-section theorem.


\section{Poisson Cross-section Theorem}

We are now ready to prove the Poisson cross-section theorem. Let $(M,\Pi)$ be a Poisson manifold with a Hamiltonian $G$-action, where $G$ is a Lie group with Lie algebra $\gfrak$. Let $\mu:M\to \gfrak^*$ be the associated equivariant moment map.

\begin{lemma} \label{lemma:kermu:poisson}
We have that $\ker\mu_*(p)\cap\im \Pi=\Pi(\gfrak_M(p)^\circ)$ and $\gfrak_p^\circ \subseteq \im\mu_*(p)$ for any $p\in M$.
\end{lemma}

\begin{Proof}
Since
\[
\mu_*\Pi(\lambda)(\xi)=d\muhat(\xi)(\Pi(\lambda))=-\lambda(\xi_M(p)), \quad \text{for $\lambda\in T_p^*M$ and $\xi\in\gfrak$,}
\]
then $\ker\mu_*(p)\cap\im \Pi=\Pi(\gfrak_M(p)^\circ)$ follows. On the other hand, if $\xi\in\gfrak$ is such that for any $v\in T_pM$ we have $\mu_*v(\xi)=d\muhat(\xi)(v)=0$, then $d\muhat(\xi)_p=0$ and $\xi_M(p)=\Pi(d\muhat(\xi))(p)=0$. By \Cref{Group-Actions:isotropy-algebra:p-tangent-orbit}, $\xi\in\gfrak_p$ and therefore $(\im\mu_*(p))^\circ\subseteq \gfrak_p$.
\end{Proof}

\begin{theorem}[Poisson Cross-section] \label{thm:poisson-cross-section}
Let $\lambda\in\gfrak^*$ and $Z\subseteq \gfrak^*$ a submanifold perfectly transverse to $G\cdot \lambda$ at $\lambda$, that is, such that $T_\lambda Z \oplus T_\lambda (G\cdot \lambda)=\gfrak^*$. Then, if $p\in M$ with $\mu(p)=\lambda$, there is a neighborhood $U$ of $p$ such that $\mu^{-1}(Z)\cap U$ is a Poisson transversal of $M$.
\end{theorem}

\begin{Proof}
The same reasoning as in the Guillemin-Sternberg cross-section theorem proves that in a neighborhood of $p$, $\mu^{-1}(Z)$ is a submanifold of $M$ such that $W_p:=T_p\mu^{-1}(Z)=\mu_*^{-1}(T_\lambda Z)$. To see that it is a Poisson transversal it suffices to see that $T_pM=W_p+\Pi(W_p^\circ)$.

Let $v\in T_pM$. By \Cref{Moment-Maps:interesting-properties} we know that
\[
\im\mu_*(p)=T_\lambda Z\cap \im\mu_*(p)\oplus T_\lambda(G\cdot \lambda),
\]
and so there are $w\in W_p$ and $\xi\in\gfrak$ such that $\mu_*v=\mu_*w+\mu_*\xi_M(p)$. Hence $v=w+\xi_M(p)+u$ with $u\in\ker\mu_*(p)$. Since $\ker\mu_*(p)\subseteq W_p$, we conclude that $T_pM=W_p+\gfrak_M(p)$. It remains to show that $\gfrak_M(p)\subseteq W_p+\Pi(W_p^\circ)$, or, equivalently, $W_p^\circ\cap \Pi^{-1}(W_p)\subseteq \gfrak_M(p)^\circ$.

As in the Guillemin-Sternberg case, since $T_\lambda Z\cap T_\lambda(G\cdot \lambda)=0$, we have that
\begin{align*}
W_p\cap \gfrak_M(p)&=\{\xi_M(p):\xi\in\gfrak, \mu_*\xi_M(p)\in T_\lambda Z\} =\{\xi_M(p):\xi\in\gfrak, \mu_*\xi_M(p)=0\}\\
& = \gfrak_M(p)\cap \ker\mu_*(p).
\end{align*}
By \Cref{lemma:kermu:poisson}, $\Pi(\gfrak_M(p)^\circ)\subseteq W_p$, and so
\[
W_p^\circ\cap \Pi^{-1}(W_p)=W_p^\circ\cap \Pi^{-1}(W_p\cap\gfrak_M(p))=W_p^\circ\cap \Pi^{-1}(\gfrak_M(p)\cap \ker\mu_*(p)).
\]
Let $\beta\in\Pi^{-1}(\gfrak_M(p)\cap \ker\mu_*(p))$. Then for any $\xi\in\gfrak$ we obtain
\[
\beta(\xi_M(p))=\beta(\Pi(d\muhat(\xi)))=-d\muhat(\xi)(\Pi(\beta))=-\mu_*\Pi(\beta)(\xi)=0.
\]
Hence $\Pi^{-1}(\gfrak_M(p)\cap \ker\mu_*(p))\subseteq \gfrak_M(p)^\circ$ and this ends the proof.
\end{Proof}

As noted before, this is just the beginning towards a Poisson implosion construction. There is yet much way to go. Indeed, let us recall what are the successive steps we have taken until reaching the conclusion that $\Mimpl$ partitions into symplectic manifolds such that its reduction by the $T$-action equals the reduction of $M$ by the $G$-action:
\begin{enumerate}
	\item The Lerman-Meinrenken-Tolman-Woodward cross-section theorem, \Cref{thm:lmtw-cross-section}, which generalizes the Guillemin-Sternberg cross-section theorem, \Cref{thm:guillemin-sternberg-cross-section}, allowed us to write $\Mimpl$ as the disjoint union of possibly singular symplectic quotients; \label{item:resumen-1}
	\item the Sjamaar-Lerman reduction theorem, \Cref{thm:sjamaar-lerman}, ensured that these symplectic quotients can be further decomposed into regular symplectic quotients, by using the decomposition of $M$ into orbit type submanifolds;
	\item symplectic reduction by stages let us see that $M\dbar_\lambda G$ and $\Mimpl\dbar_\lambda T$ are symplectomorphic for $\lambda\in\tfrak^*_+$. \label{item:resumen-3}
\end{enumerate}

The path to follow seems rather natural: Marsden-Ratiu reduction by Poisson distributions gives a natural framework for reducing a Poisson manifold by a Hamiltonian action
; hence, it seems plausible that the results just enunciated can be properly generalized to the Poisson setting, obtaining a decomposition of the Poisson imploded cross-section into Poisson manifolds such that \cref{item:resumen-3} applies (substituting ``symplectomorphic'' by ``canonically diffeomorphic''). To see \cref{item:resumen-3} in the Poisson setting, we would somehow need Poisson reduction by stages, and, happily, this has already been considered, see {\cite[Sect. 5.3]{reduction-by-stages}.

Our contribution is the first part of the generalization of step \ref{item:resumen-1}. We have shown that the Poisson manifolds into which the Poisson imploded cross-section will partition will be reduced spaces of Poisson transversals of the original manifold. This makes sense because, as has been shown, Poisson transversals inherit a natural Poisson structure from the ambient manifold.

\settocdepth{chapter}

\appendix
\chapter*{Appendix} \label{appendix}
\addcontentsline{toc}{chapter}{Appendix}
\markboth{Appendix}{Appendix}
We gather here some basic notions and conventions used throughout the text.

For any vector bundle $E\to M$ on a smooth manifold $M$ of dimension $n$ we write $\Gamma(E)$ for the vector space of sections of the bundle. For the particular case of alternating covariant and contravariant $k$-tensor bundles we use the symbols
\[
\Omega^k(M):=\Gamma(\Lambda^kT^*M)\qquad\text{and}\qquad \Xfrak^k(M):=\Gamma(\Lambda^kTM),
\]
and call them \textbf{differential $\boldsymbol{k}$-forms} and \textbf{$\boldsymbol{k}$-vector fields}, respectively. For the graded vector spaces of all differential forms and all multivector fields we use 
\[
\Omega^\bullet(M):=\bigoplus_{k=0}^n\Omega^k(M)\qquad\text{and}\qquad \Xfrak^\bullet(M):=\bigoplus_{k=0}^n\Xfrak^k(M).
\]
\bigskip

If $M$ and $N$ are manifolds and $F:M\to N$ is a smooth map, we write $F_*$ for the associated lifted map to the tangent bundles $F_*:TM\to TN$, and call it the \textbf{pushforward} of $F$. We define also the \textbf{pullback on differential forms} by $F$ to be the map $F^*:\Omega^\bullet(N)\to \Omega^\bullet(M)$ which sends $\alpha\in\Omega^k(N)$ to $F^*\alpha\in\Omega^k(M)$ given by
\[
(F^*\alpha)_p(v_1,\dots,v_k):=\alpha_{F(p)}(F_*v_1,\dots,F_*v_k), \quad\text{for every $v_1,\dots,v_k\in T_pM$.}
\]

If $F$ is diffeomorphism, we can define two more associated maps: the lifted map to cotangent bundles $F^*:T^*N\to T^*M$, which we call the \textbf{pullback} of $F$, defined by $F^*\xi=\xi\circ F_*$ for $\xi\in T^*N$, and acts fiberwise, meaning that $F^*(T_p^*N)\subseteq T_{F^{-1}(p)}^*M$; and the \textbf{pushforward on multivector fields} $F_*:\Xfrak^\bullet(M)\to\Xfrak^\bullet(N)$ which sends $X\in\Xfrak^k(M)$ to $F_*X\in\Xfrak^k(N)$ given by
\[
(F_*X)_p(\lambda_1,\dots,\lambda_k):=X_{F^{-1}(p)}(F^*\lambda_1,\dots,F^*\lambda_k), \quad\text{for every $\lambda_1,\dots,\lambda_k\in T_p^*N$.}
\]
Note that if $\alpha\in\Omega^k(N)$ and $X_1,\dots,X_k\in\Xfrak(M)$, then $F^*\alpha(X_1,\dots,X_k)=\alpha(F_*X_1,\dots, F_*X_k)\circ F$, and if $X\in\Xfrak^k(M)$ and $\alpha^1,\dots,\alpha^k\in\Omega(N)$, then $F_*X(\alpha^1,\dots,\alpha^k)=X(F^*\alpha^1,\dots,F^*\alpha^k)\circ F^{-1}$.
\bigskip

Remember that a graded algebra is an algebra $(V,\wedge)$ decomposed as a direct sum of vector spaces
\[
V=\bigoplus_{k=0}^\infty V_k
\]
satisfying $V_k\wedge V_l\subseteq V_{k+l}$. A linear map $D:V\to V$ is said to be of degree $a$ if $D(V_k)\subseteq V_{k+a}$. If, in addition, it fulfills 
\[
D(\alpha\wedge\beta)=D\alpha\wedge \beta +\epsilon^{ak}\alpha\wedge D\beta,\quad\text{ for $\alpha\in V_k$, $\beta\in V_l$ and $\epsilon=\pm 1$,}
\]
then it is called a \textbf{derivation} for $\epsilon=1$ or \textbf{antiderivation} for $\epsilon=-1$.

There are three natural derivations on $(\Omega^\bullet(M),\wedge)$: the \textbf{exterior differential} $d$, which is a degree 1 antiderivation with $d^2=0$ such that for any $f\in\Cinf(M)$, $df$ is the differential of $f$; the \textbf{inner differential} or inner product $i_X$ with $X\in\Xfrak(M)$, a degree $-1$ antiderivation with $i_X^2=0$, defined by
\[
i_X\alpha(X_2,\dots,X_k):=\alpha(X,X_2,\dots,X_k),\text{ for $\alpha\in\Omega^k(M)$ and $X_i\in\Xfrak(M)$;}
\] 
and the \textbf{Lie derivative} $\Lcal_X$, a degree 0 derivation commuting with $d$ defined by
\[
\Lcal_X\alpha:=\derev{}{t}{0}\theta_t^*\alpha,\text{ for any $\alpha\in\Omega^\bullet(M)$ and where $\theta_t$ is the flow of $X$.}\footnote{The flow is defined such that $\displaystyle X(p)=\derev{}{t}{0}\theta_t(p)$.}
\]
They are related by Cartan's magic formulas
\[
\Lcal_X=di_X+i_Xd=(d+i_X)^2\qquad\text{and}\qquad i_{[X,Y]}=\Lcal_Xi_Y-i_Y\Lcal_X.
\]
The Lie derivative can be extended in the obvious way to act on any covariant tensor field $A\in\Gamma((T^*M)^{\otimes k})$.
\bigskip

We recall that the \textbf{Lie bracket} on $\Xfrak(M)$
\[
[X,Y]:=XY-YX,\quad\text{for $X,Y\in\Xfrak(M)$,}
\]
endows $\Xfrak(M)$ with the structure of a Lie algebra. It can also be computed as the Lie derivative of $Y$ along $X$:
\[
[X,Y](p)=\Lcal_XY(p):=\derev{}{t}{0}(\theta_{-t})_*Y(\theta_t(p)),\text{ where $\theta_t$ is the flow of $X$.}
\]

For a 1-form $\alpha\in\Omega^1(M)$ and a 2-form $\omega\in\Omega^2(M)$, and vector fields $X,Y,Z\in\Xfrak(M)$, the following formulas for the exterior differential are useful
\[
d\alpha(X,Y)=X\alpha(Y)-Y\alpha(X)-\alpha([X,Y]),
\]
and
\[
d\omega(X,Y,Z)=X\omega(Y,Z)+Y\omega(Z,X)+Z\omega(X,Y)-\omega([X,Y],Z)-\omega([Y,Z],X)-\omega([Z,X],Y).
\]

The Lie bracket on vector fields can be uniquely extended to a bracket on all of $\Xfrak^\bullet(M)$, called the \textbf{Schouten-Nijenhuis bracket}, such that
\begin{enumerate}[label=(SN\arabic*),leftmargin=3\parindent]
	\item it is graded-commutative:
	\[
	[X,Y]=(-1)^{pq}[Y,X],\quad\text{for $X\in\Xfrak^p(M)$ and $Y\in\Xfrak^q(M)$;}
	\]
	\item for each $X\in\Xfrak^p(M)$, the map $\ad X:=[X,\cdot]$ is a degree $p-1$ antiderivation with respect to $\wedge$:
	\[
	[\Xfrak^p(M),\Xfrak^q(M)]\subseteq \Xfrak^{p+q-1}(M)
	\]
	and
	\[
	[X,Y\wedge Z]=[X,Y]\wedge Z+(-1)^{(p-1)q}Y\wedge [X,Z], \quad\text{for $Y\in\Xfrak^q(M)$ and $Z\in\Xfrak^r(M)$;}
	\]
	\item it satisfies the graded Jacobi identity:
	\[
	(-1)^{(p-1)r}[X, [Y, Z]] + (-1)^{(q-1)p}[Y, [Z, X]] + (-1)^{(r-1)q}[Z, [X, Y ]] = 0,
	\]
	for $X\in\Xfrak^p(M)$, $Y\in\Xfrak^q(M)$ and $Z\in\Xfrak^r(M)$;
	\item if $X,Y\in\Xfrak(M)$, $[X,Y]$ is the usual Lie bracket.
\end{enumerate}

Explicitly, it can be computed for $X_1,\dots,X_p,Y_1,\dots,Y_q\in\Xfrak(M)$ as
\begin{multline*}
[X_1\wedge\dots\wedge X_p,Y_1\wedge\dots\wedge Y_q]=\\
=(-1)^{p+1} \sum_{i=1}^p\sum_{j=1}^q (-1)^{i+j}[X_i,Y_j]\wedge X_1 \wedge \dots\wedge \hat{X}_i \wedge \dots\wedge X_p\wedge Y_1\wedge \dots\wedge\hat{Y}_j\wedge\dots\wedge Y_q.
\end{multline*}
\bigskip

Let $G$ be a Lie group. Recall that its Lie algebra $\gfrak$ is defined as the vector space $T_eG$ endowed with the following Lie bracket: $[\xi,\eta]:=[X_\xi,X_\eta](e)$, where $X_\xi\in\Xfrak(G)$ is the left-invariant vector field associated to $\xi\in\gfrak$, that is, $X_\xi(g):=L_{g*}\xi$. Because any left-invariant vector field on $G$ is complete ({\cite[Thm. 9.18]{lee}}), for any $\xi\in\gfrak$, the maximal integral curve of $X_\xi$ passing through the identity element $e\in G$ is defined for all $\R$, and we call it $\gamma_\xi$. We can then define the \textbf{exponential map} as the map
\[
\begin{array}{rrcl}
\exp: & \gfrak & \longrightarrow & G \\
& \xi & \longmapsto & \gamma_\xi(1)
\end{array}.
\]
The integral curve can hence be rewritten as $\gamma_\xi(t)=\exp(t\xi)$. 

If $\varphi:G\to H$ is a Lie group homomorphism and $\varphi_*:\gfrak\to\hfrak$ its pushforward at the identity, then $\varphi_*$ is a Lie algebra homomorphism and $\varphi\circ\exp=\exp\circ\varphi_*$ (see {\cite[Props. 4.1.5, 4.1.7]{abraham-marsden}}).

\bigskip

An $m$-dimensional \textbf{foliation} of an $n$-dimensional manifold $M$ is a partition of $M$ into disjoint connected $m$-dimensional submanifolds $\{L_a\}_a$, called leaves, such that at every point $p\in M$ there is some coordinate chart $(U,(x^i))$ in which the leaves $U\cap L_a$ are the preimages of the $m$-planes $\{(x^{m+1},\dots,x^n)=(k^{m+1},\dots, k^n)\}$, for fixed $k^{m+1},\dots, k^n\in\R$.
\bigskip

Let $P$ and $M$ be smooth manifolds, $G$ a Lie group acting smoothly and freely on $P$ on the right and $\pi:P\to M$ a smooth surjective map. Then $(P,M,G,\pi)$ is a \textbf{principal $\boldsymbol{G}$-bundle} over $M$ if for every $p\in M$ there is a neighborhood $U$ of $p$ and a diffeomorphism $\phi:\pi^{-1}(U)\to U\times G$ such that
\begin{enumerate}[label=(PB\arabic*),leftmargin=3\parindent]
	\item $\pr_1\circ \phi=\pi$, where $\pr_1:U\times G\to U$ is the projection onto the first factor, and \label{PrincBund:1}
	\item $\phi$ is $G$-equivariant with respect to the action of $G$ on $U\times G$ by right translations on the $G$ factor. \label{PrincBund:2}\qedhere
\end{enumerate}

\printbibliography[heading=bibintoc, title={References}]

\end{document}